\documentclass[11pt]{article}
\usepackage{epsfig,amsmath,latexsym,amssymb}
\usepackage{amscd,pdflscape,amsthm}
\usepackage{multirow}
\usepackage{fancyvrb, bm}
\usepackage{graphicx, subcaption, caption, listings, url}
\usepackage{enumitem, framed}
\usepackage[top=2.8cm,bottom=2.3cm,textwidth=16cm]{geometry}
\def\R{{\mathbb R}}  
\def\N{{\mathbb N}}  
\def\p{{\mathbb P}}

\def\E{{\mathbb E}}

\def\mcE{{\mathcal E}}  
\def\mcV{{\mathcal V}}

\newcommand{\Remm}[1]{}
\newtheorem{theo}{Theorem}[section]
\newtheorem{lemma}[theo]{Lemma}

\newtheorem{model ass}[theo]{Model Assumptions}

\theoremstyle{definition}
\newtheorem{remark}[theo]{Remark}

\numberwithin{equation}{section}
\renewcommand{\thesection}{\arabic{section}}

\newenvironment{Proof}[1][\unskip]{\footnotesize \textbf{Proof #1.}}{ \hfill  $\square$ \par}

  {\begin{tcolorbox}[colback=gray!15,colframe=gray!15, fonttitle=\bfseries,
 				 left=1mm, right=1mm, top=.01mm, bottom=.1mm, middle=1mm]
	\begin{lemma}}%
  {\end{lemma}\end{tcolorbox}}

  {\begin{tcolorbox}[colback=gray!15,colframe=gray!15, fonttitle=\bfseries,
 				 left=1mm, right=1mm, top=.01mm, bottom=.1mm, middle=1mm]
	\begin{theo}}%
  {\end{theo}\end{tcolorbox}}

\begin{document}
\author{Philippe Deprez\footnote{RiskLab, Department of Mathematics, 
ETH Zurich, 8092 Zurich, Switzerland} \qquad 
Mario V.~W\"uthrich$^\ast$\footnote{Swiss Finance Institute SFI Professor}}
\date{\today}
\title{Construction of Directed Assortative Configuration Graphs}
\maketitle

\begin{abstract} 
Constructions of directed configuration graphs 
based on a given bi-degree distribution were introduced  
in random graph theory some years ago. 
These constructions lead to graphs  
where the degrees of two nodes belonging to the same 
edge are independent. 
However, it is observed that many real-life networks 
are assortative, meaning that edges tend to connect low degree nodes 
with high degree nodes, or variations thereof. 
In this article we provide an explicit algorithm to construct 
directed assortative configuration graphs based on a given 
bi-degree distribution and an arbitrary pre-specified assortativity.  
\end{abstract}


\section{Introduction}\label{Section: Introduction}
Random graphs are used to model large networks 
that consist of particles, called nodes, which are possibly 
linked  to each other by edges. 
The study of random graphs goes back to the works of~\cite{Erdos} and~\cite{Gilbert}. 
Since then, numerous random graph models 
have been introduced and studied in the literature. 
For an overview we refer the reader to~\cite{Bollobas, Durrett, Remco5,Remco4}. 
Empirical studies of large data sets of real-life networks have shown that in many 
cases the degrees of two nodes belonging to the same edge 
are not independent (where the degree of a node is 
defined to be the number of edges attached to it).  
It is observed that in some types of real-life networks the degree of a node 
is positively related to the degrees of its linked neighbors, 
while in other situations the degree of a node 
is negatively related to the degrees of its linked neighbors. 
This property is called {\it assortativity} or {\it assortative mixing}. 
It has been discovered by~\cite{Bech, Cont, Soramaki} that financial 
networks typically show negative assortativity and that the 
strength of the assortativity influences 
the vulnerability of the financial network to shocks, see also~\cite{Hurd}. 
In contrast, social networks tend to be 
positive assortative, see for instance~\cite{Newman2}. 
More examples of assortative networks are presented in~\cite{Newman} 
and~\cite{Remco2}, where also quantities to measure the assortativity 
in networks are proposed.  
On the other hand, there is only little literature on explicit constructions 
of random graphs showing assortative mixing. 
For example,~\cite{Stanton} and~\cite{Bassler} study the construction 
of graphs based on a given graphical degree sequence, 
and~\cite{Bloznelis} analyzes assortativity in random intersection graphs. 

\begin{figure}
\begin{center}
\includegraphics[width=7cm]{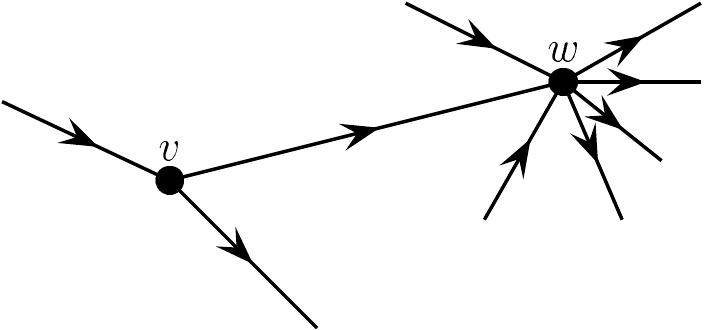}
\end{center}
\caption{\footnotesize
Node $v$ is of type $(1,2)$ and node $w$ is of type $(3,4)$. 
Edge $e=\langle v,w \rangle$ is of type $(2,3)$.} 
\label{Figure: Node Type}
\end{figure}

\begin{figure}
\centering
\begin{subfigure}{.33\textwidth}
  \centering
  \includegraphics[width=.9\linewidth]{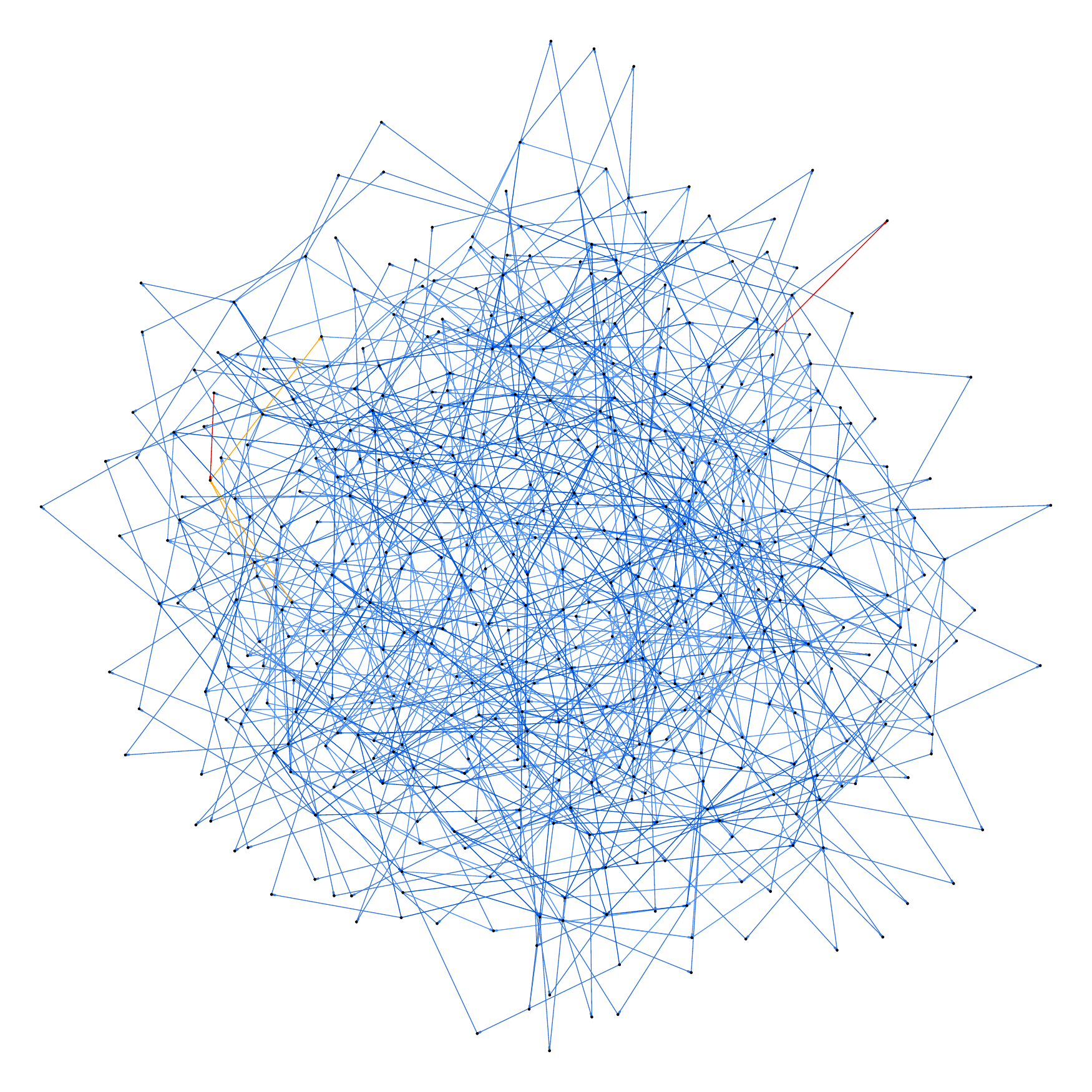}
  \caption{\footnotesize negative assortativity}
\end{subfigure}%
\begin{subfigure}{.33\textwidth}
  \centering
  \includegraphics[width=.9\linewidth]{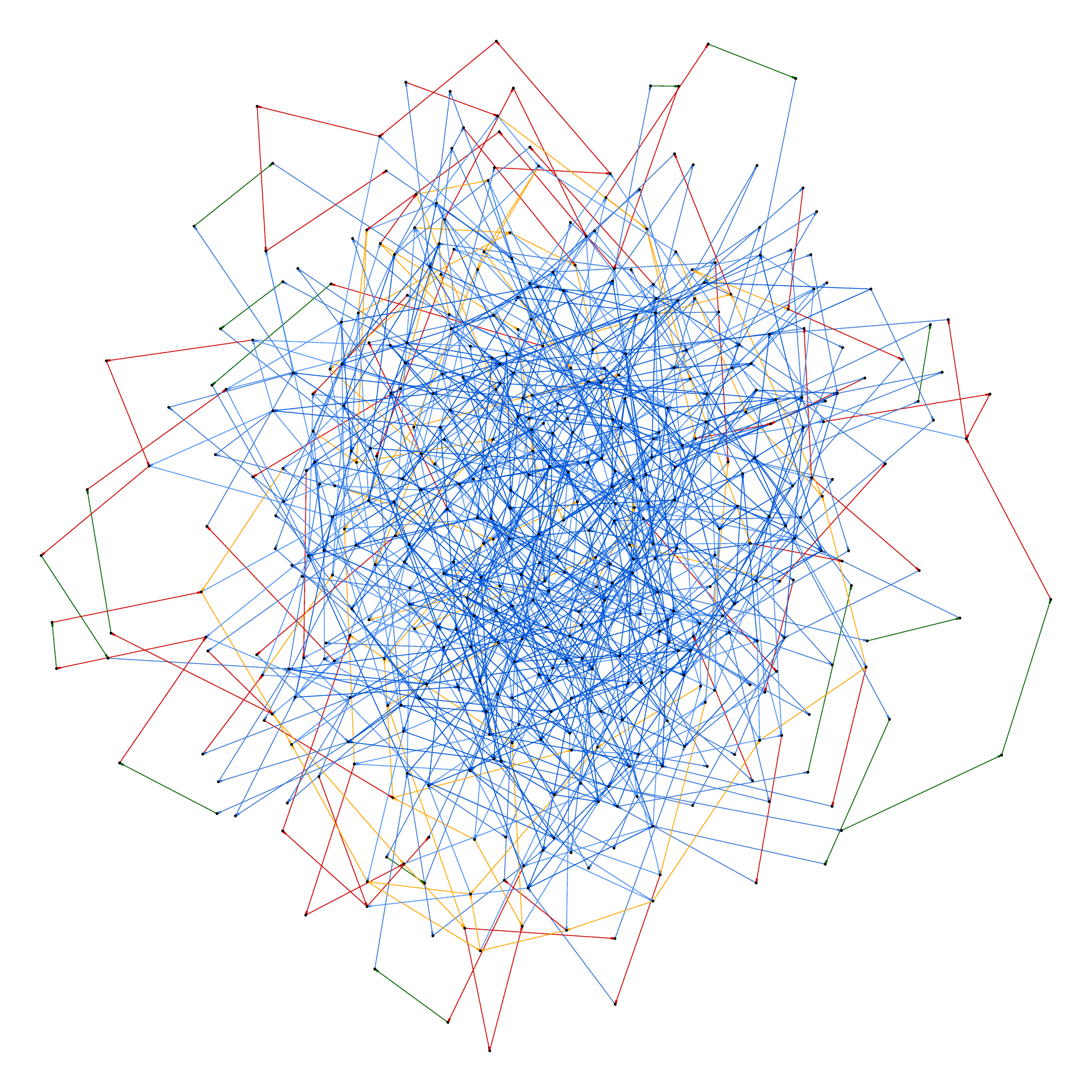}
  \caption{\footnotesize non-assortativity}
\end{subfigure}%
\begin{subfigure}{.33\textwidth}
  \centering
  \includegraphics[width=.9\linewidth]{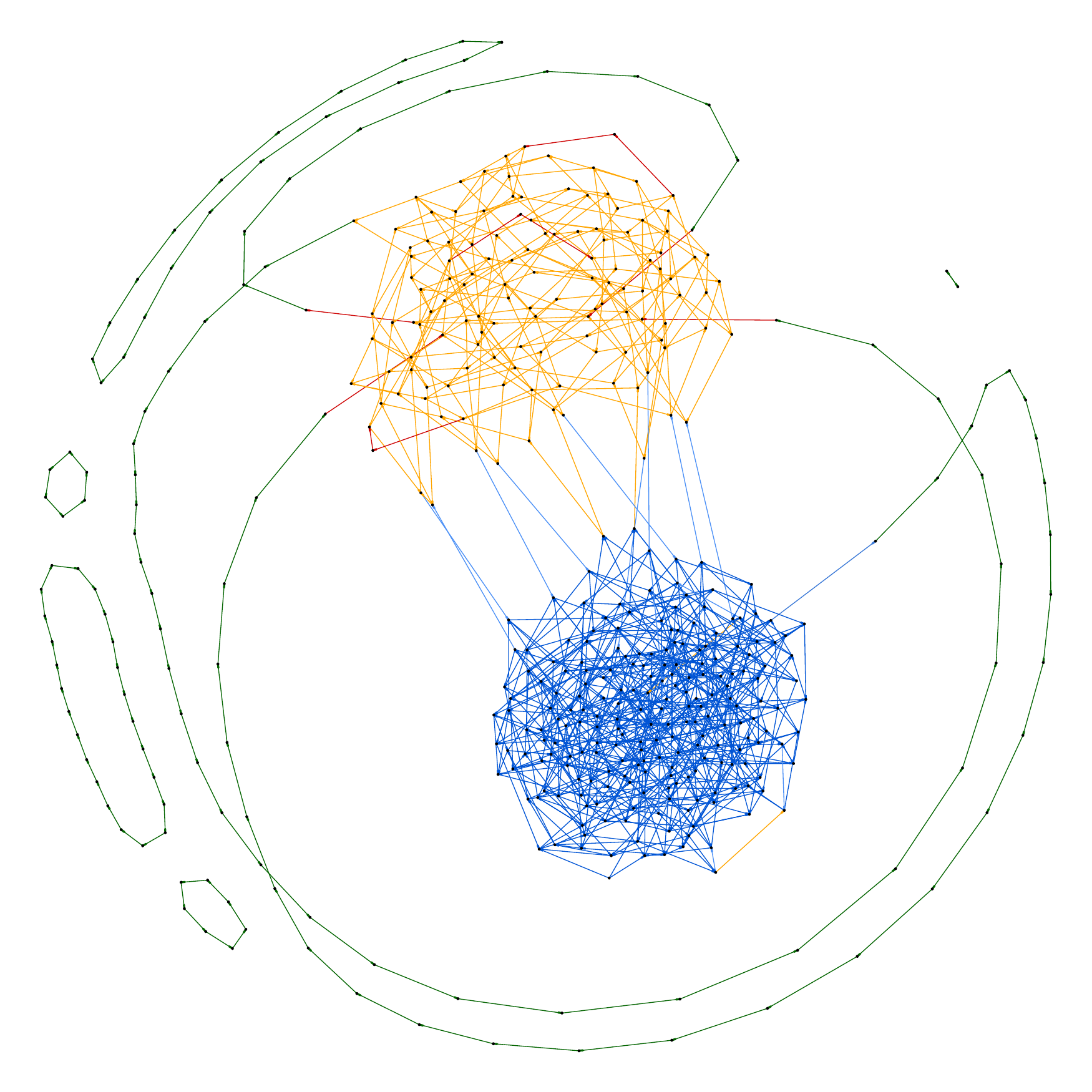}
  \caption{\footnotesize positive assortativity}
\end{subfigure}
\caption{\footnotesize
Three graphs generated by our algorithm 
given in Section~\ref{Section: Algorithm} 
with  $N=500$ nodes. 
They all have the same 
node-type distribution $P$ but different 
edge-type distributions $Q$. 
There are three different types of nodes present in each graph: 
$(1,1)$, $(2,2)$ and $(3,3)$. 
Edges of identical type are colored the same. 
Edges that are arriving or leaving a highest-degree node 
are colored in different shades of blue. 
These edges are mainly present in the negative assortative case (a). 
In the positive assortative case (c), mainly nodes of the same type are connected.  
These edges are colored dark blue, orange and green. 
All other possible edges are colored red, which significantly appear only in (b).
}
\label{Figure: Introduction}
\end{figure}

However, established constructions of directed random graphs 
based on a 
given bi-degree distribution, called {\it configuration graphs}, 
lead to {\it non-assortative} graphs,  
see for instance the construction presented in \cite{Chen}.
Here, the bi-degree of a node $v$ 
is a tuple $(j_v,k_v)$, where $j_v$ is the number of edges arriving  
at node $v$ (called {\it in-degree}) and $k_v$ is the number of edges 
leaving from node $v$ (called {\it out-degree}), 
and we say that node $v$ is of type $(j_v,k_v)$, 
see Figure~\ref{Figure: Node Type} for an illustration. 
In this article we extend the non-assortative 
construction presented in~\cite{Chen} by  
giving an explicit algorithm 
which allows to construct directed configuration graphs 
with a pre-specified assortativity based on 
a concept introduced in~\cite{Gleeson}. 
Namely, \cite{Gleeson} proposed to specify 
the graph not only through their node-types, 
but also through their edge-types. 
We define the type of an edge $e=\langle v,w \rangle$ connecting 
node $v$ to node $w$  by a tuple  
$(k_e,j_e)$ with $k_e$ denoting the out-degree of node $v$ and 
$j_e$ denoting the in-degree of node $w$, 
see Figure~\ref{Figure: Node Type} for an illustration. 
This notion of edge-types is directly related to 
the notion of assortativity. 
In the positive assortative case $k_e$ is positively 
related to $j_e$ meaning that edges tend to connect 
nodes having similar degrees, and accordingly for the 
negative assortative case. 
If $k_e$ is independent of $j_e$, then the graph is 
non-assortative. 
This motivated \cite{Gleeson}  to 
construct  directed assortative configuration graphs 
based on a given node-type distribution $P$ 
describing the nodes {\it and} from a given edge-type distribution $Q$ 
describing the edges, while different choices of  $Q$ result in different 
types of assortativity in the constructed graphs, 
see Figure~\ref{Figure: Introduction} for examples. 
Nevertheless, there was not an explicit construction given in~\cite{Gleeson}.  
The aim of this article is to construct random graphs where 
nodes and edges follow pre-specified given bivariate distributions $P$ and $Q$, 
and we give a precise mathematical meaning to the choice of these distributions. 

~

Let us first interpret the meaning of the distributions $P$ and $Q$ in more detail.  
Node-type distribution $P$ has the following interpretation. 
Assume we have a large directed network and we choose 
{\it at random} a node $v$ of that network, then the type $(j_v,k_v)$  
of $v$ has distribution $P$.   
Similarly, edge-type distribution $Q$ should be understood as follows. 
If we choose {\it at random} an edge $e$ of a large network, 
then its type $(k_e,j_e)$ has distribution $Q$.  
This concept of $P$ and $Q$ distributions seems straightforward, 
however, it needs quite some care in order to 
give a rigorous mathematical meaning to these distributions, the difficulty 
lying in the ``randomly'' chosen node and edge obeying $P$ and $Q$, 
respectively:  
the graph as total induces dependencies between nodes and edges 
which implies that the exact distributions can only be obtained 
in an asymptotic sense (this will be seen in the construction below).  

~

We give an explicit algorithm to construct a directed 
assortative configuration graph with a given number of nodes and 
based on  distributions $P$ and  $Q$,   
and we prove that the type of a randomly chosen node of the resulting graph 
converges in distribution to $P$ 
as the size of the graph tends to infinity. 
Similarly, the type of a randomly chosen edge converges
in distribution to $Q$. 
These convergence results give a rigorous mathematical meaning to $P$ and $Q$ 
in line with their interpretation given above. 
The proposed algorithm allows for self-loops and multiple edges. 
In order to obtain a simple graph we delete all self-loops and 
multiple  edges, and we show that the convergence results 
still hold true for the resulting simple graph. 
Recently, an alternative approach to construct assortative 
configuration graphs based on given distributions $P$ and $Q$ 
was proposed in~\cite{Hurd2}, see also \cite{Hurd}, 
using techniques from~\cite{Wormald}. 
Our construction is different from~\cite{Hurd2} 
and more in the spirit of~\cite{Bollobas2, Remco3, Chen}.  
Moreover, we give a rigorous mathematical meaning 
to the given distributions $P$ and $Q$ which 
relies on the law of large numbers only. 

~

In Section~\ref{Section: Model} we introduce the model 
and state our main results. Section~\ref{Section: Algorithm} 
specifies the algorithm to 
generate directed assortative configuration graphs. 
The implementation of the algorithm 
in the programming language \textsf{R} 
can be downloaded from: 
\begin{equation}\label{Link}
\text{\small\url{https://people.math.ethz.ch/~wueth/Papers/AssortativeConfigurationGraphs.R}}
\end{equation}
In Section~\ref{Section: Examples} we illustrate examples 
of  assortative configuration graphs generated by our algorithm 
showing different assortative mixing, and we study their 
empirical assortativity coefficients as well as their 
empirical node- and edge-type distributions. 
The proofs of the results are given in Section~\ref{Section: Proofs}.


\section{Model and main results}\label{Section: Model}
Consider fixed finite integers $J\ge 1$ and $K \ge 1$ 
which describe the maximal in- and out-degree 
of a node, respectively. 
For  $l\ge0$ and $n\ge l$ define $[n]_l=\{l,\ldots,n\}$. 
For $j \in [J]_0$ and $k \in [K]_0$ 
we say that node $v$ is of type $(j,k)$ if the in-degree of $v$ is $j$ 
and the out-degree of $v$ is $k$. 
For $k \in [K]_1$ and $j \in [J]_1$ 
we say that a directed edge $e=\langle v,w \rangle$ 
is of type $(k,j)$ if the out-degree of $v$ is $k$ and the in-degree 
of $w$ is $j$. 
Figure~\ref{Figure: Node Type} illustrates the notions  
of node- and edge-types. 
In the remainder, letter $j$ always 
refers to in-degree and letter $k$ to out-degree. 
Consider two bivariate probability distributions 
\begin{eqnarray*}
P&=&(p_{j,k})_{j\in[J]_0,k\in[K]_0} 
\qquad \text{with $\sum_{j\in[J]_0,k\in[K]_0}p_{j,k}=1$;}\\
Q&=&(q_{k,j})_{k\in[K]_1, j\in[J]_1} 
\qquad \text{with $\sum_{k\in[K]_1,j\in[J]_1}q_{k,j}=1$}.
\end{eqnarray*}
We call $P$ node-type distribution 
and $Q$ edge-type distribution. 
We denote the marginal distributions of $P$ and $Q$ respectively by 
\begin{eqnarray*}
p_j^- &=& \sum_{k'\in[K]_0} p_{j,k'}
\qquad \text{and} \qquad 
p_k^+ ~=~ \sum_{j'\in[J]_0} p_{j',k},
\qquad \text{$j\in[J]_0$ and $k\in[K]_0$; } \\
q_k^+ &=& \sum_{j'\in[J]_1} q_{k,j'} 
\qquad \text{and} \qquad
q_j^- ~=~ \sum_{k'\in[K]_1} q_{k',j},
\qquad \text{$k\in[K]_1$ and $j\in[J]_1$}. 
\end{eqnarray*}
In the remainder, superscript ``$-$'' always 
refers to in-degree and superscript ``$+$''  to out-degree. 
For instance, $(p_j^-)_{j\in[J]_0}$ denotes the in-degree distribution of nodes.   
Observe that in a given graph  
the number of edges $e=\langle v,w \rangle$ with out-degree of $v$ being $k\in[K]_1$ 
is equal to $k$ times the number of nodes having out-degree $k$, 
and similarly for the number of nodes having in-degree $j\in[J]_1$. 
This relation between nodes and edges implies that 
we cannot choose $P$ and $Q$ independently of each other 
to achieve that nodes and edges in the constructed graph 
follow $P$ and $Q$, respectively. 
We therefore assume that $P$ and $Q$ satisfy the following consistency conditions, 
see also~\cite{Gleeson} and~\cite{Hurd2}, 
which implies that the above observation holds true 
in expectation in graphs where nodes and edges follow distributions $P$ and $Q$,
respectively. 
\begin{align}
\label{Condition 2}\tag{C1}
q_k^+&~=~kp_k^+/z, \qquad  k\in [K]_1;
\\\label{Condition 3}\tag{C2}
q_j^-&~=~jp_j^-/z,    \qquad j\in [J]_1,
\end{align}
with {\it mean degree} $z=\sum_{k \in [K]_0} k p_k^+$. 
Observe  that conditions~\eqref{Condition 2} and~\eqref{Condition 3}  
require that $z=\sum_{k \in [K]_0} k p_k^+ =\sum_{j \in [J]_0} j p_j^->0$. 
This says that, in expectation, the sum of in-degrees 
equals the sum of out-degrees if nodes and edges 
follow  distributions $P$ and $Q$, respectively.

\begin{remark}
We assume uniformly bounded degrees. 
A generalization to unbounded degrees is possible but not straightforward: 
conditions~\eqref{Condition 2} and~\eqref{Condition 3} need to be 
fulfilled and the rate of decay of $p_k^+$ and $p_j^-$ needs 
a careful specification so that the results below still remain true. 
\end{remark}

~

\begin{remark}\label{Remark: variation}
We have defined assortativity through edge-types, which, 
for an edge $e=\langle v,w \rangle$ connecting 
node $v$ to node $w$, is defined to be the tuple  
$(k_e,j_e)$ with $k_e$ denoting the out-degree of  $v$ and 
$j_e$ denoting the in-degree of  $w$. 
There are three other possibilities to define the type of an edge $e$, 
for instance, by the tuple $(k_e,k'_e)$ with $k_e$ and $k'_e$ denoting the out-degree of  $v$ and 
$w$, respectively. 
We comment on these variations of edge-types in more detail in Appendix~\ref{Appendix}. 
\end{remark}

~

Given the number of nodes $N\in\N$ and given distributions $P$ and $Q$ 
satisfying~\eqref{Condition 2} and~\eqref{Condition 3}  
the goal is to construct a graph such that the following 
statement is true in an asymptotic sense as the size $N$ of 
the graph tends to infinity:  
the type of a randomly chosen node has distribution $P$ and 
the type of a randomly chosen edge has distribution $Q$. 
The following theorem shows that this is indeed the case 
for graphs constructed by the algorithm provided in Section~\ref{Section: Algorithm} 
and, hence, the theorem gives an explicit mathematical meaning to $P$ and $Q$. 

~

\begin{theo}\label{Lemma: random degrees}
Fix $s\in\N$. 
Let $(j_{v_1},k_{v_1}),\ldots,(j_{v_s},k_{v_s})$ be the 
types of $s$ randomly chosen nodes of the graph 
generated by the algorithm provided in Section~\ref{Section: Algorithm}.  
Then,  
\begin{equation*}
\big((j_{v_1},k_{v_1}),\ldots,(j_{v_s},k_{v_s})\big) \stackrel{d}{\longrightarrow} 
\big((j'_1,k'_{1}),\ldots,(j'_s,k'_{s})\big),
\qquad \text{as $N\to\infty$}, 
\end{equation*}
where $(j'_1,k'_{1}),\ldots,(j'_s,k'_{s})$ are $s$ independent  
random variables having distribution~$P$. 
Similarly, the types of $s$ randomly chosen edges converge in 
distribution, as $N\to\infty$, to a sequence of $s$ independent 
random variables having distribution~$Q$. 
\end{theo}

~

If we consider a graph where nodes and edges 
follow distributions $P$ and $Q$, respectively,  
then we expect that the relative number of nodes of type $(j,k)$ 
is close to $p_{j,k}$ 
and that the relative number of edges of type $(k,j)$ is close 
to $q_{k,j}$. 
Theorem~\ref{Lemma: empirical distributions 2} below makes this statement precise 
for graphs constructed by the algorithm provided 
in Section~\ref{Section: Algorithm}. 
To formulate the theorem, denote by $\mcV_{j,k}$ 
the number of nodes of type $(j,k)$, $j \in [J]_0$ and $k\in [K]_0$, 
and by $\mcE_{k,j}$ the number of edges 
of type $(k,j)$, $k \in [K]_1$ and $j\in [J]_1$, of the constructed graph of size $N$. 
The total number of edges is denoted by $\mcE$. 
Theorem~\ref{Lemma: empirical distributions 2} says that 
the relative frequencies $\mcV_{j,k}/N$ and 
$\mcE_{k,j}/\mcE$ converge to $p_{j,k}$ and $q_{k,j}$, 
respectively, in probability as $N\to\infty$. 

~

\begin{theo}\label{Lemma: empirical distributions 2}
For the random graph constructed by the algorithm provided 
in Section~\ref{Section: Algorithm} 
we have for any $\varepsilon>0$ 
\begin{equation*}
\lim_{N\to\infty}
\p\left[
\sum_{j \in [J]_0, k\in [K]_0} \left|
\frac{\mcV_{j,k}}{N}
- p_{j,k}  
\right| 
+
\sum_{k \in [K]_1, j\in [J]_1} \left|
\frac{\mcE_{k,j} }{\mcE}
- q_{k,j}   
\right|
> \varepsilon
\right]
~=~0. 
\end{equation*}
\end{theo}

~

The algorithm provided in Section~\ref{Section: Algorithm} 
generates a graph possibly not being simple, i.e.~it 
may contain self-loops and multiple edges. 
To obtain a simple graph we delete (erase) all self-loops and multiple edges,  
and we call the resulting graph {\it erased configuration graph}. 
The following theorem states that the asymptotic results 
still hold true for the erased configuration graph. 

~

\begin{theo}\label{Theorem: simple}
The results of Theorem~\ref{Lemma: random degrees}
and Theorem~\ref{Lemma: empirical distributions 2}   
still hold true for the erased configuration graph, 
based on the algorithm provided in Section~\ref{Section: Algorithm}. 
\end{theo}


\section{Construction of directed assortative configuration graphs}\label{Section: Algorithm}
The algorithm to construct directed assortative configuration graphs 
starts from the work of \cite{Chen}, where the authors construct a directed random 
graph with $N\in\N$ nodes based on given in-degree and 
given out-degree distributions in the following way. 
They assign to each node independently an in-degree and 
an out-degree according to the given distributions, also 
independently for different nodes. 
Some degrees are then modified if the sum of in-degrees 
differs from the sum of out-degrees so that these sums of degrees are equal, and 
the sample is only accepted if the number of 
modifications is not too large. 
Finally, in-degrees are randomly paired with out-degrees. 
Note that this construction leads to a {\it non-assortative} configuration graph 
and the in- and out-degree of a given node are independent. 
The construction of an assortative configuration graph 
is more delicate  
since in-degrees cannot be randomly paired with out-degrees. 
In our construction we generate node-types  
using directly  node-type distribution $P$. 
Independently of the node-types we generate $zN$ 
edges having independent 
edge-types according to distribution $Q$. 
Finally, we match  in- and out-degrees of nodes 
with  edges of corresponding types. 
In general, the matching cannot be 
done exactly, but with high probability the number of types 
that need to be changed accordingly 
is small for large $N$,  
due to consistency conditions~\eqref{Condition 2} and~\eqref{Condition 3}. 
We first describe the algorithm in detail and then 
comment on each step of the algorithm below. 

~

\noindent
{\bf Algorithm to construct directed assortative configuration graphs.}\nopagebreak

\noindent
Assume maximal degrees $J,K\ge1$ and  
two probability distributions 
$P$ and $Q$ satisfying~\eqref{Condition 2} and~\eqref{Condition 3} 
with mean degree $z$ are given. 
Choose $\delta \in(1/2,1)$ fixed. 
Choose $N\in\N$ so large that there 
exists  $N'\in\N$ with  
$N=N' + 2\left\lceil N^\delta \right \rceil + \max\{J^2,K^2\}$,  
and set $N''=N' + \left\lceil N^\delta\right \rceil$. 
Here, $\left\lceil x \right \rceil$ denotes the smallest 
integer larger than or equal to $x\in\R$. 
\begin{enumerate}[label={}, leftmargin=10pt, rightmargin=10pt]
\item 
{\bf Step 1.}  
Assign to each node $v=1,\ldots,N'$ independently a node-type $(j_v,k_v)$ 
according to distribution $P$. 
Generate edges $e=1,\ldots, \lceil zN'' \rceil$ 
having independent edge-types $(k_e,j_e)$ according to distribution $Q$,  
independently of the node-types. 
Define
\begin{eqnarray*}
n_k^+ &=& \sum_{v=1}^{N'} 1_ {\{ k_v=k \}}
\quad \text{ and } \quad 
e_k^+ ~=~ \left\lceil \frac{1}{k}\sum_{e=1}^{\lceil zN'' \rceil} 1_ {\{ k_e=k \}}\right\rceil
\quad \text{ for all $k\in[K]_1$};
\\
n_j^- &=& \sum_{v=1}^{N'} 1_ {\{ j_v=j \}}
\quad \text{ and } \quad 
e_j^- ~=~ \left\lceil \frac{1}{j}\sum_{e=1}^{\lceil zN'' \rceil} 1_ {\{ j_e=j \}}\right\rceil
\quad\text{ for all $j\in[J]_1$}.  
\end{eqnarray*}
Let $A_N$ be the event on which we have 
\begin{eqnarray*}
\left|  n_k^+ - p_k^+ N' \right| \le p_k^+N^\delta/2
&\quad \text{ and } \quad &
\left|  e_k^+ - p_k^+ N'' \right| \le p_k^+N^\delta/2
\quad \text{ for all $k\in[K]_1$}; 
\\
\left|  n_j^- - p_j^- N' \right| \le p_j^-N^\delta/2
&\quad \text{ and } \quad &
\left|  e_j^- - p_j^- N'' \right| \le p_j^-N^\delta/2
\quad\text{ for all $j\in[J]_1$}.  
\end{eqnarray*}
Proceed to Step~2 if event $A_N$ occurs. 
Otherwise, proceed to Step~5.

\item 
{\bf Step 2.}
For each $k\in[K]_1$ and each $j\in[J]_1$ do the following. 
\begin{itemize}
\item
Add $r_k^+ = ke_k^+- \sum_{e=1}^{\lceil zN'' \rceil} 1_ {\{ k_e=k \}}$ edges of type $(k,1)$; 
\item
Add $r_j^- = je_j^-- \sum_{e=1}^{\lceil zN'' \rceil} 1_ {\{ j_e=j \}}$ edges of type $(1,j)$.
\end{itemize}
Set $r^+ = \sum_{k\in[K]_1}r_k^+$ and $r^- = \sum_{j\in[J]_1}r_j^-$. 

\item 
{\bf Step 3.}
Set the type of each node in $\{N'+1,\ldots,N\}$ to $(0,0)$. 
For each $k\in[K]_1$ and each $j\in[J]_1$ do the following. 
\begin{itemize}
\item
Take the first $e_k^+ - n_k^++r^-1_{\{k=1\}}$ nodes 
in $\{N'+1,\ldots,N\}$ having out-degree $0$ and change their 
out-degrees to $k$; 
\item
Take the first $e_j^- - n_j^-+r^+1_{\{j=1\}}$ nodes 
in $\{N'+1,\ldots,N\}$ having in-degree $0$ and change their 
in-degrees to $j$.
\end{itemize}

\item 
{\bf Step 4.}
For each $k\in[K]_1$ and each $j\in[J]_1$ do the following. 
\begin{itemize}
\item
Assign to each node 
having out-degree $k$ exactly $k$ uniformly chosen  
edges $e$ of type $(k_e,j_e)$ with $k_e=k$; 
\item
Assign to every node 
having in-degree $j$ exactly $j$ uniformly chosen 
edges $e$ of type $(k_e,j_e)$ with $j_e=j$.
\end{itemize}
Proceed to Step~6. 

\item 
{\bf Step 5.}
Define node-types $(j_1,k_1)=(0,1)$, $(j_2,k_2)=(1,0)$ and  
$(j_v,k_v)=(0,0)$ for all $v=3,\ldots, N$.  
Insert an edge $e=\langle 1,2 \rangle$ that connects 
node $1$ to node $2$.  

\item 
{\bf Step 6.}
Return the constructed graph.
\end{enumerate}

~

\noindent
{\bf Explanation of the algorithm.}  \nopagebreak

\noindent
We say that node $v$ is a {\it $k$-node} if its out-degree is $k\in[K]_1$,  
and similarly we say that edge $e=\langle v,w \rangle$ is a 
{\it $k$-edge} if $v$ is a $k$-node, $k\in[K]_1$. 

\noindent
{\bf Step 1.} 
We generate only $N'$ node-types and 
we keep $2\lceil N^\delta\rceil + \max\{J^2,K^2\}$ nodes undetermined 
for possible modifications in later steps. 
The expected number of generated $k$-nodes is  
 $p_k^+ N'$ and the expected number of 
generated $k$-edges is $q_k^+\lceil zN'' \rceil$.  
Using condition~\eqref{Condition 2}, the expected number of $k$-nodes 
needed for the generated $k$-edges is 
therefore $q_k^+\lceil zN'' \rceil/k \approx p_k^+N'' > p_k^+ N'$. 
Henceforth, if $e_k^+$ is close to its 
expectation $q_k^+\lceil zN'' \rceil/k \approx p_k^+N''$, 
it dominates the number of generated $k$-nodes $n_k^+$ which is 
of order $p_k^+N' < p_k^+ N''$. 
Step~3 is then used to correct for this imbalance in 
a deterministic way, and event $A_N$ guarantees that this correction  
is possible. For receiving an efficient algorithm we would like event $A_N$ 
to occur sufficiently likely, which is exactly stated in the next lemma. 

~

\begin{lemma}\label{Lemma: acceptance set}
We have $\p\left[A_N\right] \to1$ as $N\to\infty$. 
\end{lemma}

~

\begin{figure}
\begin{center}
\includegraphics[width=12cm]{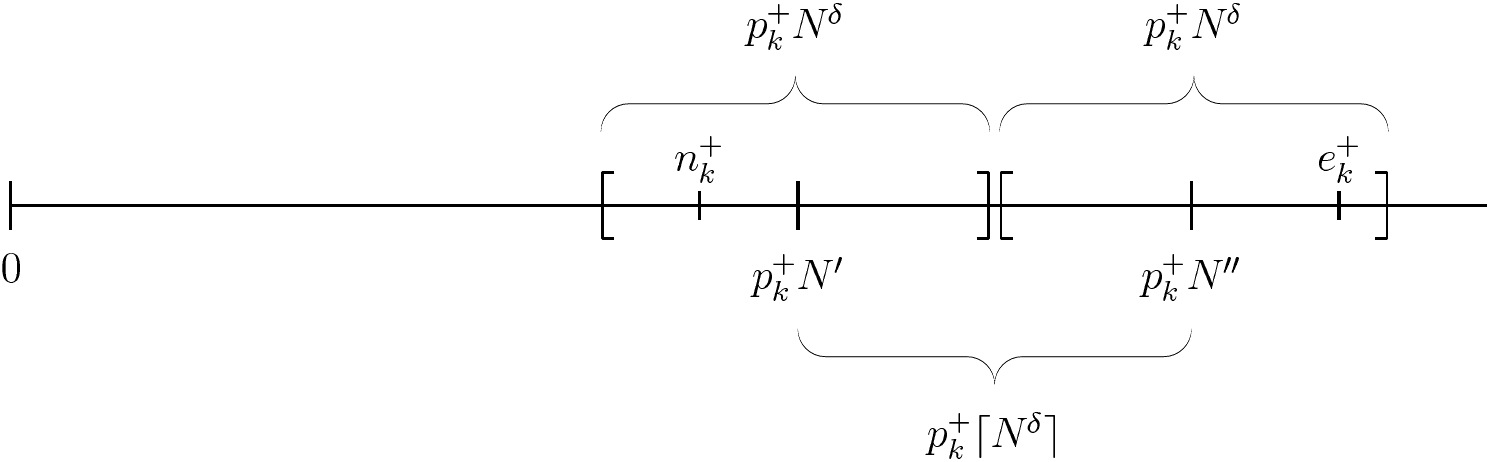}
\end{center}
\caption{\footnotesize
On event $A_N$, the number of generated $k$-nodes, $n_k^+$, 
lies in the interval of length $p_k^+N^\delta$ around $p_k^+N'$.  
The number of $k$-nodes needed for the generated $k$-edges, $e_k^+$,  
lies in the interval of length $p_k^+N^\delta$ around $p_k^+N''$. 
By definition of $N''$ the gap between the two intervals 
is of size between $0$ and $p_k^+$. 
From this it follows that there are at most $2p_k^+\lceil N^\delta \rceil $ additional 
$k$-nodes needed in order to attach all generated 
$k$-edges to $k$-nodes. 
} 
\label{Figure: Construction}
\end{figure}

Integer $N''$ is chosen in such a way that, on event $A_N$, 
the number of $k$-nodes needed 
for the generated $k$-edges dominates 
the number of generated $k$-nodes, and  
their difference 
is at most $2p_k^+\lceil N^\delta \rceil$ for each $k\in[K]_1$, 
see also Figure~\ref{Figure: Construction} for an illustration.  
Therefore, on event $A_N$, the total number of additional nodes  needed 
having a positive out-degree is 
\begin{eqnarray*}
\sum_{k\in[K]_1} \left( e_k^+ - n_k^+\right)
&\le&
2\lceil N^\delta \rceil \sum_{k\in[K]_1} p_k^+
~\le~
2\lceil N^\delta \rceil. 
\end{eqnarray*}
Hence, we have sufficiently many undetermined nodes 
in $\{N'+1,\ldots,N\}$ to which we can assign 
out-degrees accordingly in Step~3, and similarly for the in-degrees.

\noindent
{\bf Step 2.} 
In general, the number of generated $k$-edges is 
not a multiple of $k$. 
Therefore, we use Step~2 to correct for this cardinality 
by defining $r_k^+$ additional edges of type $(k,1)$. 
Note that each such edge requires a node 
having in-degree $1$. 
Therefore, in total $r^+ \le K^2$ nodes having in-degree $1$ 
are additionally needed, and similarly for the added edges of type $(1,j)$, 
$j\in[J]_1$. 
The undetermined $\max\{J^2,K^2\}$ nodes are 
exactly used to correct for the  
corresponding node-types in Step~3. 

\noindent
{\bf Steps 3 and 4.} 
We assign node-types to the undetermined nodes $N'+1,\ldots,N$  
in such a way that the total number of $k$-nodes is equal to 
the number of $k$-nodes needed for the $k$-edges 
generated in Steps~1 and~2, for each $k\in[K]_1$, and similarly 
for the total number of nodes having in-degree $j\in[J]_1$.  
Then, all cardinalities for $j$ and $k$ match and 
all edges can be randomly connected to corresponding nodes.   
After doing so, each node has 
the correct number of arriving and leaving edges according to its type. 
Note that this step allows for self-loops and multiple edges. 

\noindent
{\bf Step 5.} 
If event $A_N$ does not occur in Step~2, we 
just define a deterministic graph having one edge so that all terms in 
Theorem~\ref{Lemma: random degrees} and 
Theorem~\ref{Lemma: empirical distributions 2} 
are well-defined. 
Due to Lemma~\ref{Lemma: acceptance set} the influence of this 
deterministic graph is negligible.


\section{Discussion and examples}\label{Section: Examples}

Given a non-degenerate node-type distribution $P$ with mean degree  
$z>0$ given by $z=\sum_{k \in [K]_0} k p_k^+=\sum_{j \in [J]_0} j p_j^-$, 
we aim to find possible edge-type distributions $Q$ such that 
$P$ and $Q$ satisfy~\eqref{Condition 2} and~\eqref{Condition 3}. 
Conditions~\eqref{Condition 2} and~\eqref{Condition 3} imply that  
the marginal distributions of $Q$ are fully described by the marginal  
distributions of $P$, and 
their respective cumulative distribution functions are given by 
\begin{equation*}
Q^+(k) ~=~ \sum_{k'=1}^{k}q_{k'}^+~=~ \frac{1}{z}\sum_{k'=1}^{k}k'p_{k'}^+
\qquad \text{ and } \qquad
Q^-(j) ~=~ \sum_{j'=1}^{j}q_{j'}^-~=~\frac{1}{z}\sum_{j'=1}^{j}j'p_{j'}^-, 
\end{equation*}
for $k\in[K]_1$ and $j\in[J]_1$. 
The possible joint distributions $Q=(q_{k,j})_{k,j}$ are therefore 
given by   
\begin{eqnarray}
q_{k,j}&=& \label{Equation: C}
C\left( Q^+(k),Q^-(j) \right)
+C\left( Q^+(k-1),Q^-(j-1)  \right)
\\&&\hspace{1.5cm}-C\left( Q^+(k),Q^-(j-1)  \right)-C\left( Q^+(k-1),Q^-(j)  \right),\notag
\end{eqnarray}
where $C:[0,1]^2 \to [0,1]$ is a $2$-dimensional copula, see for instance~\cite{QRM}. 
To measure assortativity of a graph,~\cite{Newman} introduced 
the {\it assortativity coefficient} of $Q$ given by 
\begin{equation*}
\rho_Q ~=~
\frac{\sum_{k\in[K]_1}\sum_{j\in[J]_1} kj \left( q_{k,j} - q_k^+q_j^- \right)}
{\sqrt{\sum_{k\in[K]_1} k^2q_k^+ - \left( \sum_{k\in[K]_1} kq_k^+ \right)^2 }
\sqrt{\sum_{j\in[J]_1} j^2q_j^- - \left( \sum_{j\in[J]_1} jq_j^- \right)^2 }}
~\in[-1,1],
\end{equation*}
which is Pearson's correlation coefficient of distribution $Q$, 
see also~\cite{Pim} for an analysis of different types of correlations in a graph. 
By Hoeffding's identity and using  representation~\eqref{Equation: C}, 
$\rho_Q$ can be rewritten as 
\begin{equation*}
\rho_Q ~=~
\frac{
\sum_{k\in[K]_1}\sum_{j\in[J]_1} 
\big(  C\left( Q^+(k),Q^-(j) \right) - Q^+(k)Q^-(j) \big)
}
{\sqrt{\sum_{k\in[K]_1} k^2q_k^+ - \left( \sum_{k\in[K]_1} kq_k^+ \right)^2 }
\sqrt{\sum_{j\in[J]_1} j^2q_j^- - \left( \sum_{j\in[J]_1} jq_j^- \right)^2 }}.
\end{equation*}
Observe that $\rho_Q$ is determined by $P$ and $C$. 
Define the copulas 
\begin{eqnarray*}
W(u_1,u_2) &=& \max\{u_1+u_2-1,0\};
\\
M(u_1,u_2) &=& \min\{u_1,u_2\};
\\
\Pi(u_1,u_2) &=& u_1u_2,
\end{eqnarray*}
for $u_1,u_2 \in [0,1]$. 
Then, $C=W$ corresponds to the 
minimal possible assortativity coefficient $\rho_Q^- \in [-1,0]$, 
and $C=M$ corresponds to the 
maximal possible assortativity coefficient $\rho_Q^+ \in [0,1]$. 
Copula $C=\Pi$ leads to non-assortativity 
and in this case we have $q_{k,j}=kjp_k^+p_j^-/z^2$ 
for all $k\in[K]_1$ and $j\in[J]_1$. 
Note that for given $P$, $\rho_Q$ does not uniquely determine $C$. 
On the other hand, one can always find $\lambda \in [0,1]$ 
such that $\lambda W + (1-\lambda)M$ leads to a given  
assortativity coefficient $\rho_Q \in [\rho_Q^-,\rho_Q^+]$. 
This allows to construct directed assortative 
configuration graphs having any given  
assortativity coefficient that is possible for given 
node-type distribution $P$. 

~

To illustrate assortativity in an example  
we consider maximal in- and out-degree $J=K=4$ 
and 
node-type distribution $P_p=\left(p^p_{j,k}\right)$, $p\in(0,1)$, 
given by 
\begin{equation*}
P_p = 
 \begin{pmatrix}
  	0 & 	0 & 	0 		& 	0 & 	0 \\
  	0 & 	0 & 	0 		& 	0 & 	0 \\
 	0 & 	0 & 	p 		& 	0 & 	0 \\
 	0 & 	0 & 	0 		& 	0 & 	0 \\
 	0 & 	0 & 	0 		& 	0 & 	1-p
 \end{pmatrix}.
\end{equation*}
Distribution $P_p$ only allows for nodes of types $(2,2)$ and $(4,4)$, 
with respective probabilities $p$ and $1-p$, which results in a mean degree 
of $z=4-2p$. 
Clearly, these nodes can only be connected through edges 
of types $(2,2)$, $(2,4)$, $(4,2)$ and $(4,4)$. 
Since $P_p$ is diagonal, consistency conditions~\eqref{Condition 2} and~\eqref{Condition 3} 
fully specify the edge-type distribution $Q_q$ which is given by 
\begin{equation*}
Q_q = \frac{1}{2-p}
 \begin{pmatrix}
  	0 & 	0 		& 	0 & 	0 \\
 	0 & 	3p+q-2		& 	0 & 	2-2p-q \\
 	0 & 	0 		& 	0 & 	0 \\
 	0 & 	2-2p-q 		& 	0 & 	q
 \end{pmatrix},
\end{equation*}
for $q=q(p)\in [\max\{2-3p,0\},2-2p]$. 
For fixed $p\in(0,1)$, different values of $q$ lead to different 
assortativity coefficients $\rho_q=\rho_{Q_q}$. 
A straightforward calculation gives 
\begin{equation*}
\rho_q
~=~
\frac{q(2-p)-4(1-p)^2}{2p(1-p)},
\quad \text { or equivalently  } \quad 
q
~=~
\frac{2(1-p)(2 +  p(\rho_q-2))}{2-p}.
\end{equation*}
For any $p\in(0,1)$, the optimal bounds on $\rho_q$ are given by 
\begin{equation*}
- \min\left\{
\frac{p}{2(1-p)},\frac{2(1-p)}{p}
\right\}
~=~
\rho_p^-
~\le~ 
\rho_q 
~\le~
\rho_p^+
~=~
1.
\end{equation*}
Observe that $\rho_p^-=-1$ if and only if  $p=2/3$.  

\begin{figure}
\centering
\begin{subfigure}{.32\textwidth}
  \centering
  \includegraphics[width=1\linewidth]{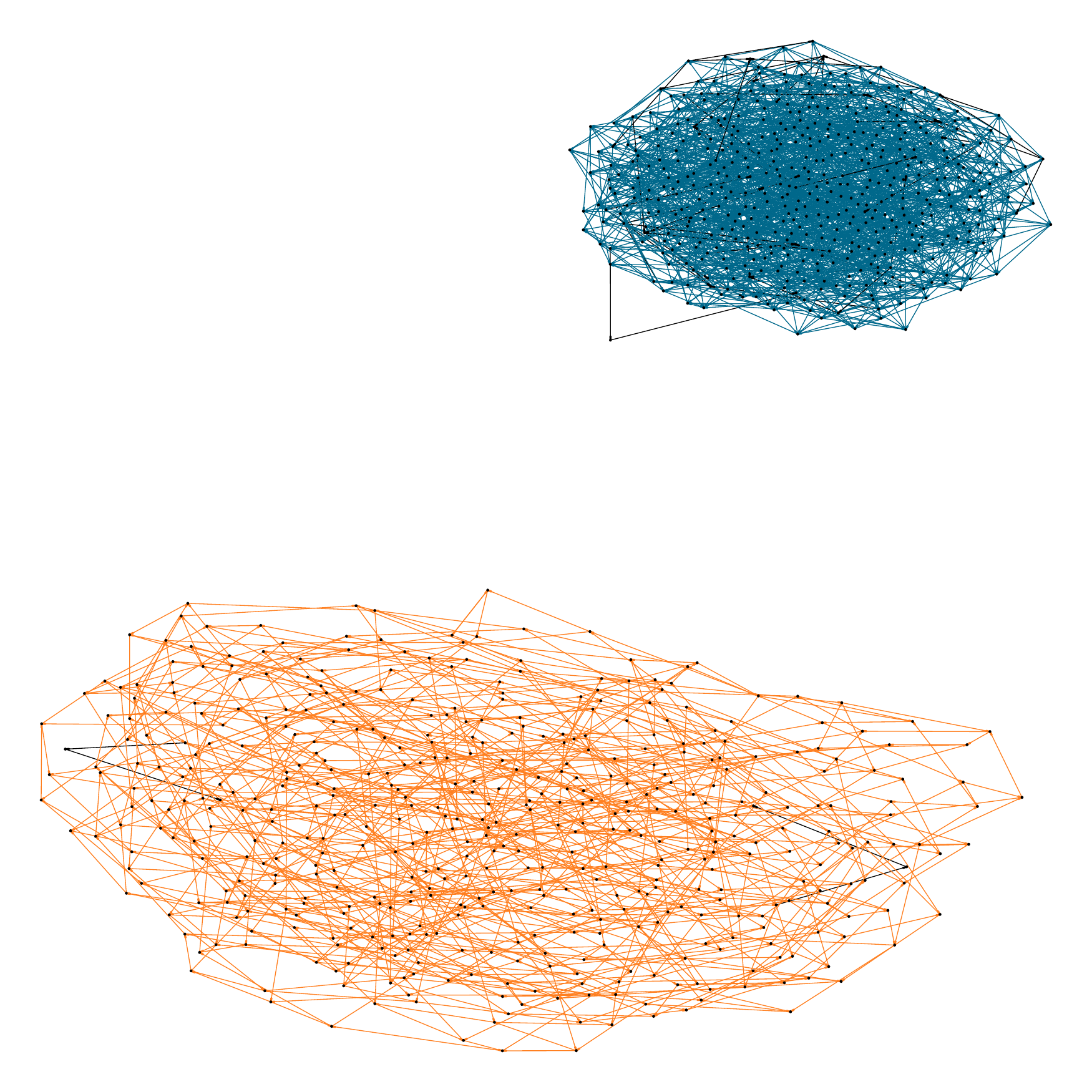}
  \caption{\footnotesize $q=1$, $\rho_{q}=1$, $\hat \rho_q=0.99$.}
\end{subfigure}%
\begin{subfigure}{.32\textwidth}
  \centering
  \includegraphics[angle=180, width=1\linewidth]{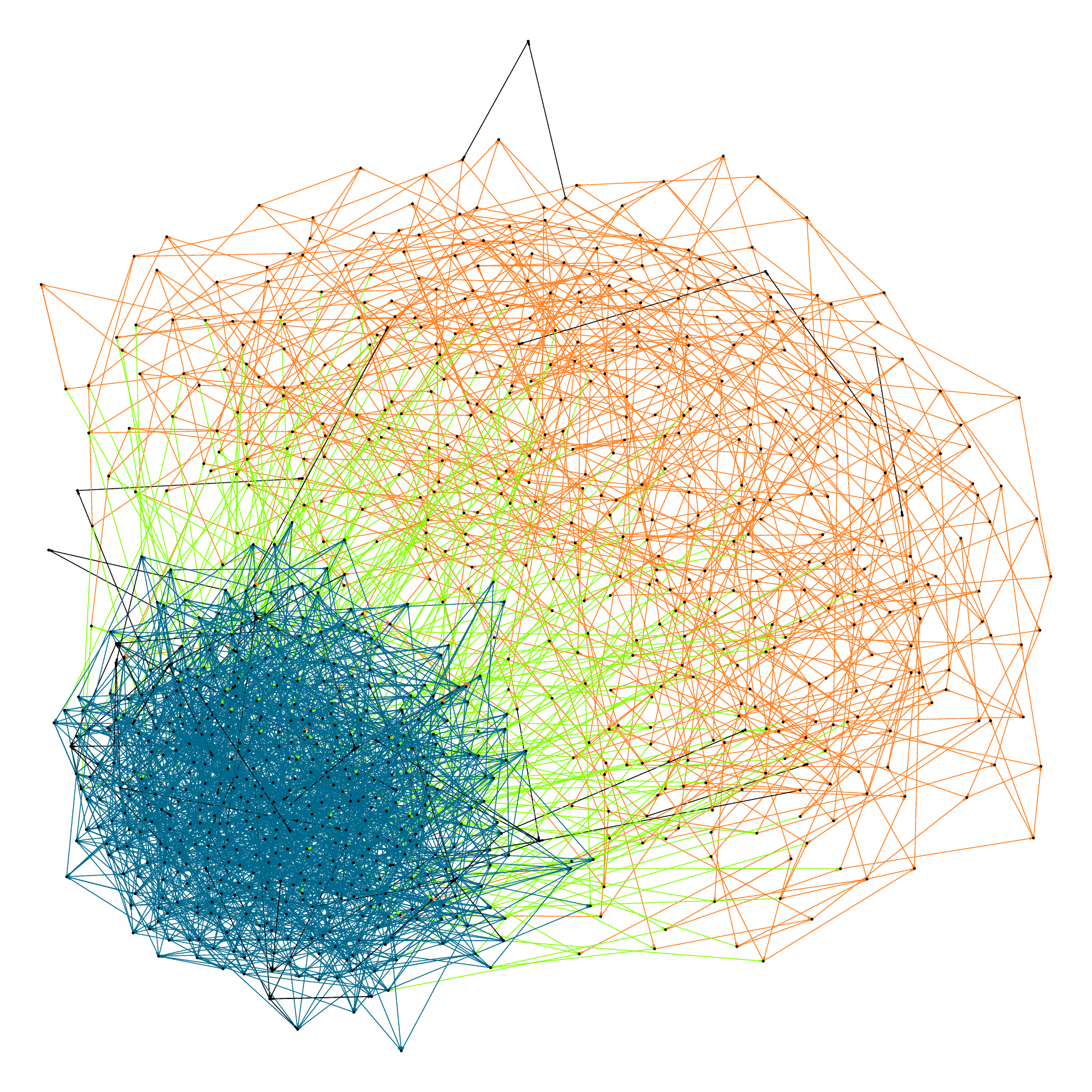}
  \caption{\footnotesize $q=\frac{14}{15}$, $\rho_{q}=0.8$, $\hat \rho_q=0.79$.}
\end{subfigure}
\begin{subfigure}{.32\textwidth}
  \centering
  \includegraphics[angle=270, width=1\linewidth]{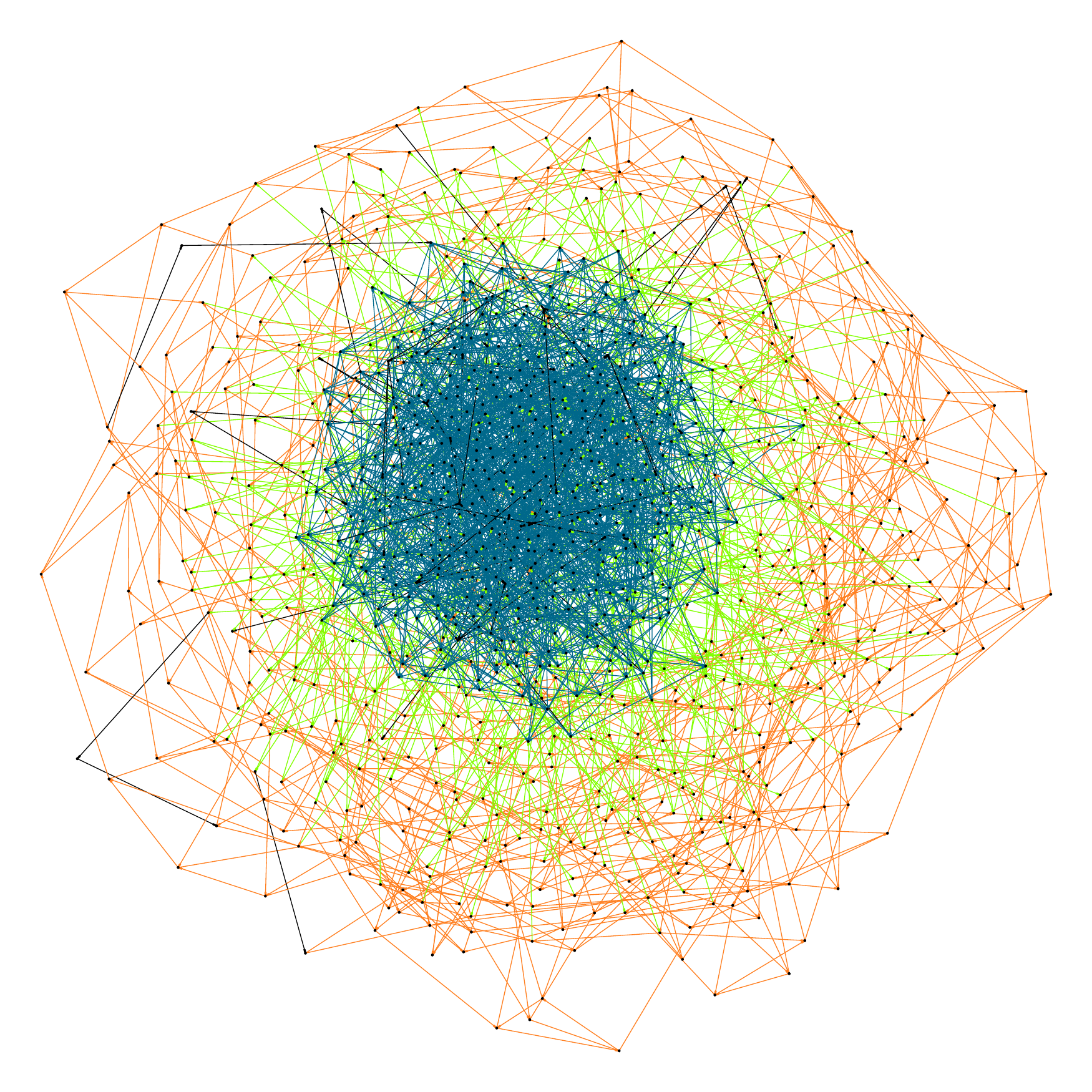}
  \caption{\footnotesize $q=\frac{13}{15}$, $\rho_{q}=0.6$, $\hat \rho_q=0.6$.}
\end{subfigure}\\
\vspace{1cm}
\begin{subfigure}{.32\textwidth}
  \centering
  \includegraphics[angle=270,width=1\linewidth]{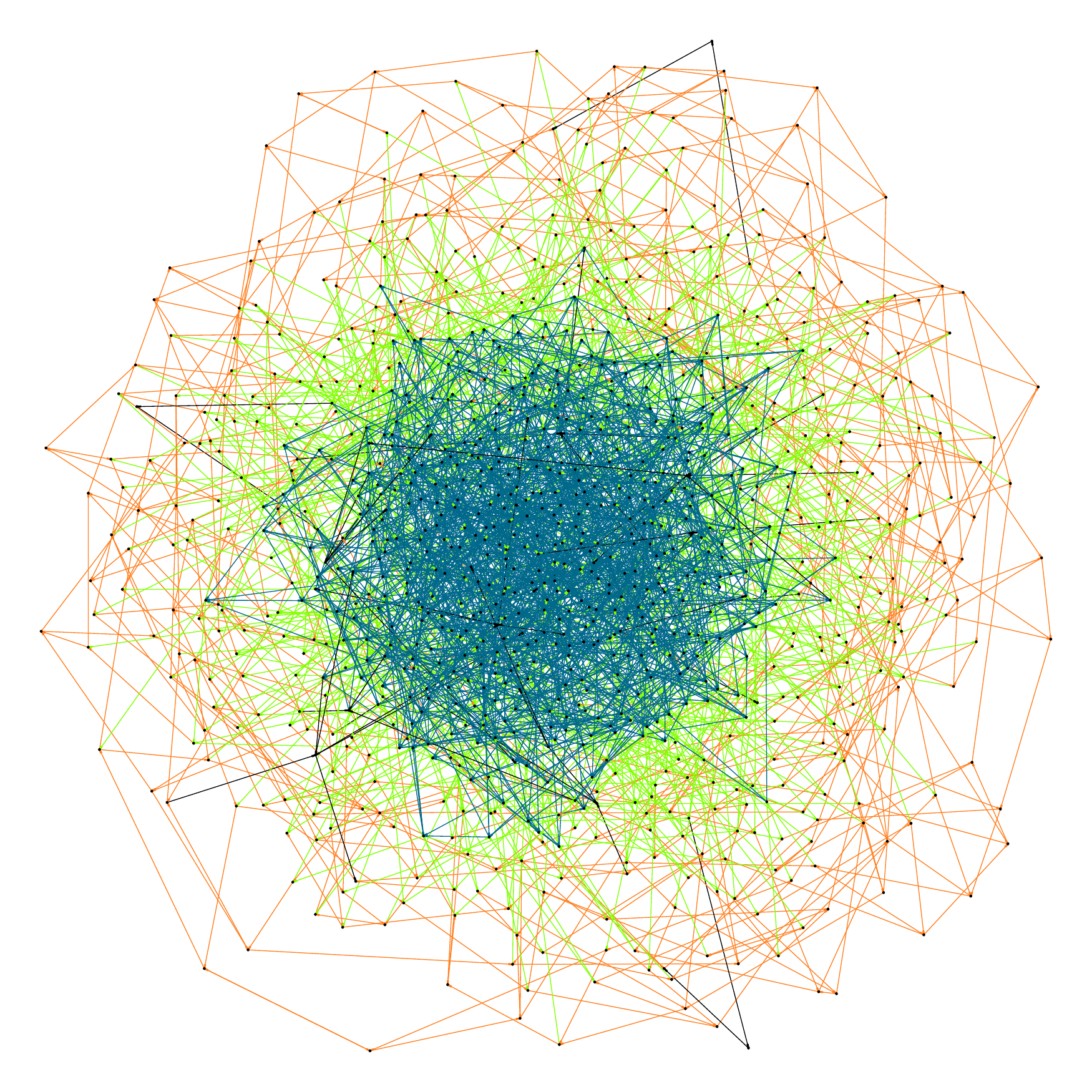}
  \caption{\footnotesize $q=0.8$, $\rho_{q}=0.4$, $\hat \rho_q=0.39$.}
\end{subfigure}%
\begin{subfigure}{.32\textwidth}
  \centering
  \includegraphics[width=1\linewidth]{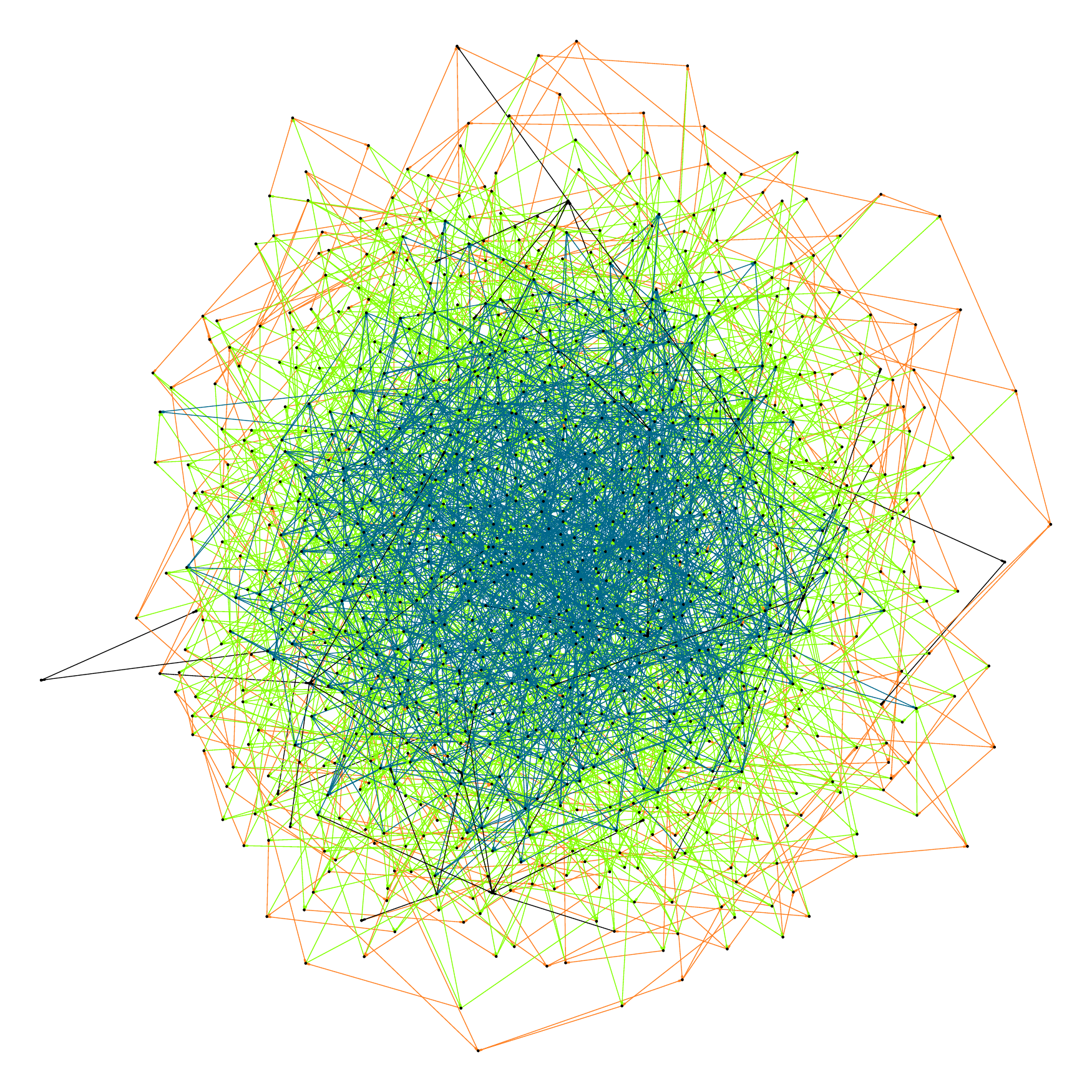}
  \caption{\footnotesize $q=\frac{2}{3}$, $\rho_{q}=0$, $\hat \rho_q=0.01$.}
\end{subfigure}
\begin{subfigure}{.32\textwidth}
  \centering
  \includegraphics[width=1\linewidth]{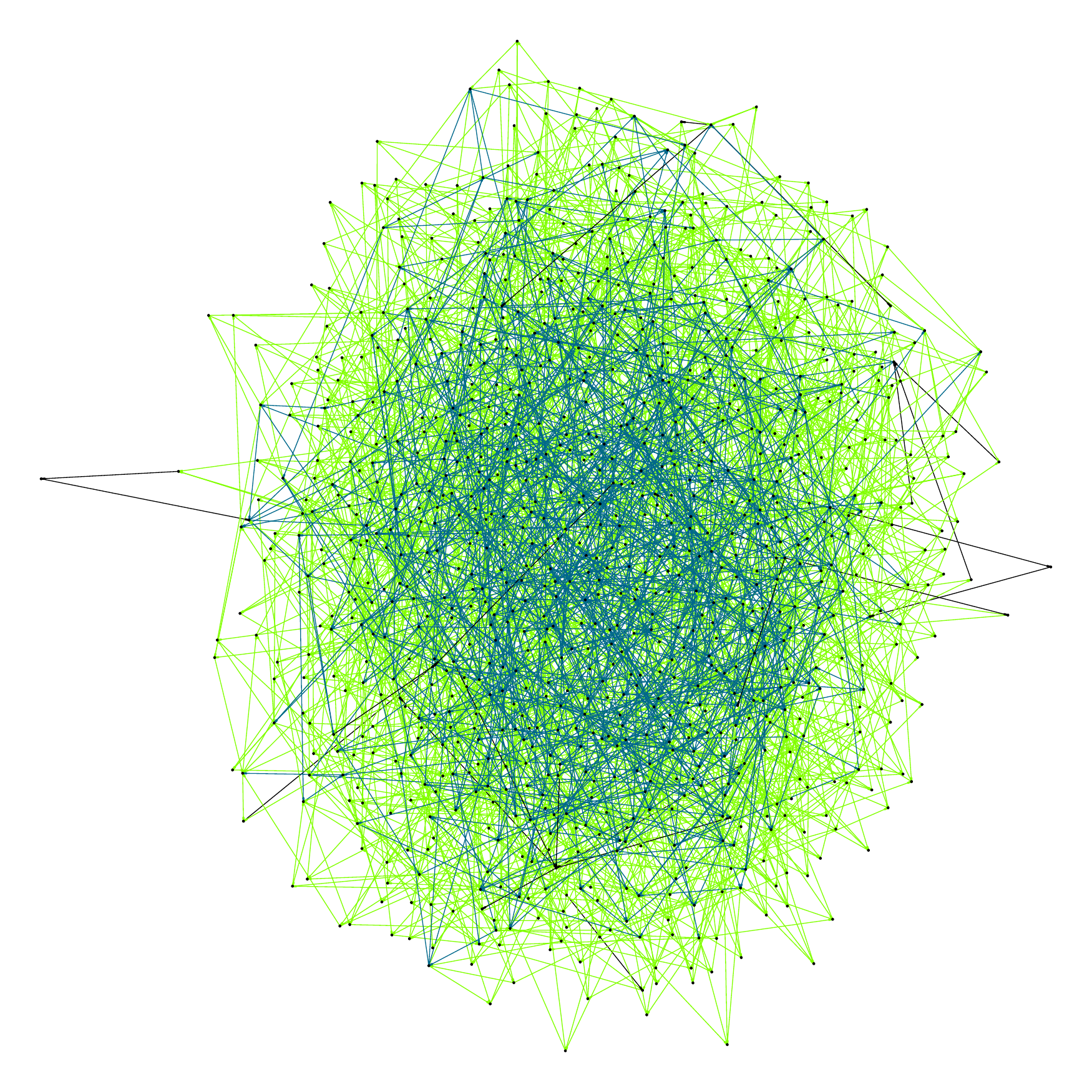}
  \caption{\footnotesize $q=0.5$, $\rho_{q}=-0.5$, $\hat \rho_q=-0.49$.}
\end{subfigure}
\caption{\footnotesize 
All graphs were generated by our algorithm in Section~\ref{Section: Algorithm} 
with $N=1000$,  
with the same node-type distribution 
$P_{0.5}={\rm diag}(0,0,0.5,0,0.5)$ but different edge-type distribution $Q_q$. 
Edges of type $(2,2)$ are colored orange, while edges 
of type $(4,4)$ are colored blue. 
Edges of types $(2,4)$ and $(4,2)$ are colored green. 
All other edges are colored black. 
}
\label{Figure: Diagonal}
\end{figure}

From now on we fix $p=0.5$, meaning that 
there are nodes of types $(2,2)$ and $(4,4)$ 
with equal probability.  
For any $q \in [0.5,1]$ or $\rho_q \in [-0.5,1]$, this leads to 
\begin{equation*}
\rho_q
~=~
3q-2 \in [-0.5,1]
\quad \text { or equivalently  } \quad 
q
~=~
\frac{\rho_q + 2}{3} \in [0.5,1].
\end{equation*}
A value of $q=1$ results in a graph with maximal 
assortativity coefficient $\rho_{1} =\rho^+_{0.5}= 1$. 
In this case, there are only edges that connect 
nodes having identical types.   
By decreasing $q$ we allow also for edges of types 
$(2,4)$ and $(4,2)$, while we reduce the probability of having 
edges of types $(2,2)$ and  $(4,4)$. 
If we decrease $q$ to its minimal value $q=0.5$, 
edges of type $(2,2)$ finally disappear and there 
are only edges of types $(2,4)$, $(4,2)$ and $(4,4)$. 
This means that if $q=0.5$, each edge is leaving from a node 
with maximal possible out-degree or is arriving at 
a node with maximal possible in-degree. 
In this case, the assortativity coefficient is negative and 
given by $\rho_{0.5} = \rho_{0.5}^-=-0.5$. 
Non-assortativity is given for $q=2/3$. 
To illustrate these different types of assortativity, 
Figure~\ref{Figure: Diagonal} shows six graphs generated 
by the algorithm given in Section~\ref{Section: Algorithm} 
with $N=1000$ nodes,  
with node-type distribution $P=P_{0.5}$ and edge-type 
distribution $Q_q$ for values of $q$ such 
that  $\rho_q\in\{-0.5,0,0.4,0.6,0.8,1\}$. 
Note that the algorithm produces node-types 
different from $(2,2)$ and $(4,4)$ due to the 
modifications on nodes $\{N'+1,\ldots,N\}$ with 
$N-N' \ll N$ as $N\to \infty$.   
For better illustration we erase self-loops and multiple edges, 
and we do not show nodes of type $(0,0)$. 
To illustrate the differences between the six generated graphs 
we color edges of identical types the same as follows: 
\begin{eqnarray*}
&&\text{edges of type $(2,2)$ are orange,}\\
&&\text{edges of types $(2,4)$ and $(4,2)$ are green,}\\
&&\text{edges of type $(4,4)$ are blue}, 
\end{eqnarray*}
and all other edges are colored black (which may arise 
by the construction and the erasure procedure). 

We analyze 
the resulting empirical distributions $\hat P=(\hat p_{j,k})=(\mcV_{j,k}/N)$ 
and $\hat Q_q=(\hat q_{k,j})=(\mcE_{k,j}/\mcE)$ for all six graphs 
in Figure~\ref{Figure: Diagonal}. 
Let us first consider the graphs without having erased 
self-loops and multiple edges.
For all simulations we have set $\delta=0.5+0.0001$ implying that 
$N-N'=80$ nodes out of the total $1000$  do not have a bi-degree 
generated from $P$. 
This  implies that the sum of the components of the difference 
$\hat P - \text{diag}(0,0,\hat p_{2,2},0,\hat p_{4,4})$ 
is at most $0.08$ for all generated graphs. 
For all the graphs in Figure~\ref{Figure: Diagonal} we observe empirically that the 
sum is at most $0.075$. 
Moreover, the differences $|\hat p_{2,2}- p_{2,2}|$ and $|\hat p_{4,4}- p_{4,4}|$ 
are both at most $0.043$ for all six graphs. 
For $\hat Q_q$ we have  $|\hat q_{k,j}- q_{k,j}| \le 0.02$ for 
all possible edge-types $(k,j)$ and all six graphs. 
From this we conclude that already for a comparably 
small graph of $N=1000$ nodes we obtain very accurate results 
(note that the results are exact for $N\to\infty$). 

\begin{figure}[t]
\begin{center}
\includegraphics[width=10cm]{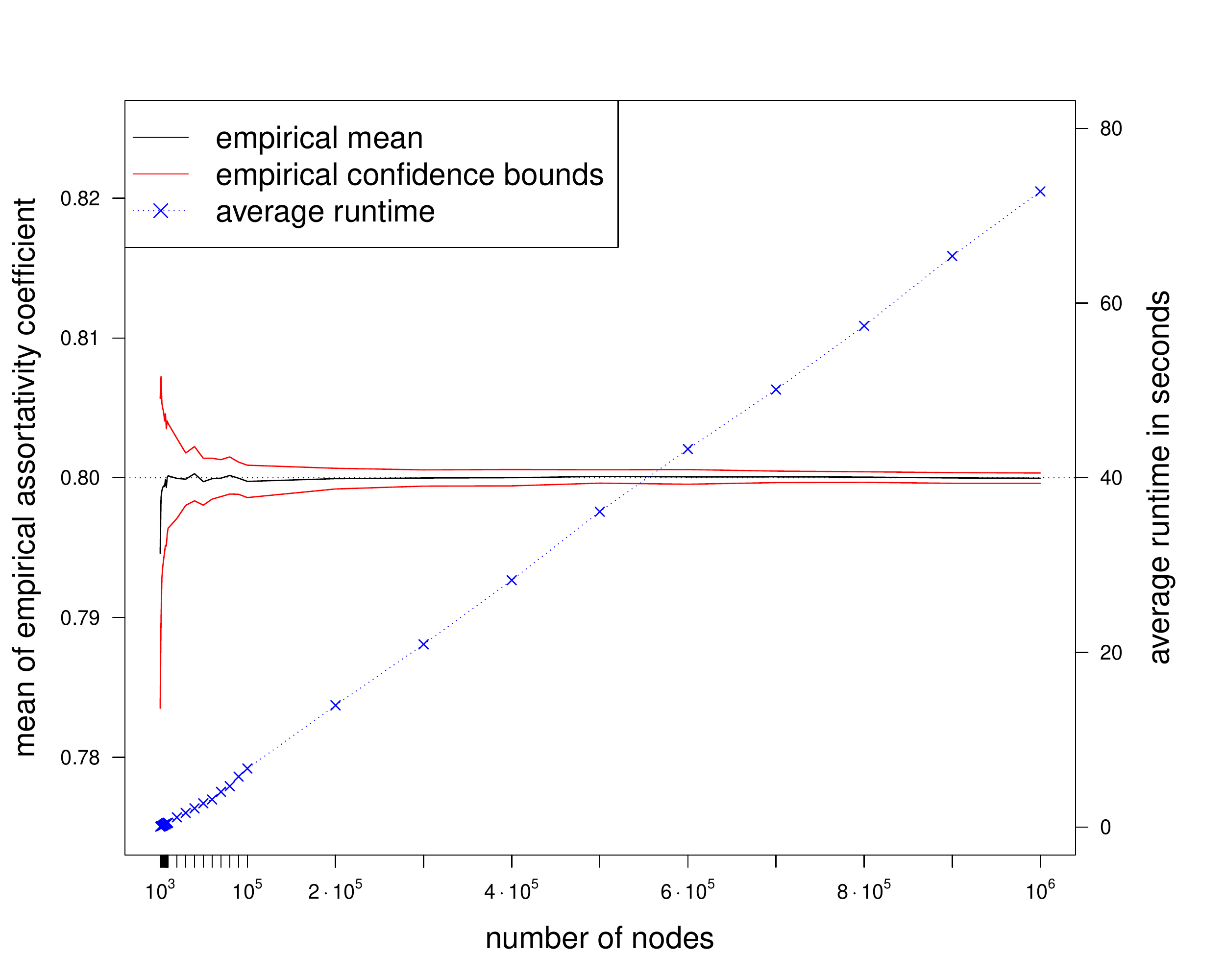}
\end{center}
\caption{\footnotesize
Mean of empirical assortativity coefficient with confidence bounds of 
one empirical standard deviation for different values of $N$ and 
$100$ simulations for each $N$ (left axis). 
Each cross indicates the average runtime in seconds for the simulation of
a graph for fixed $N$ (right axis). 
} 
\label{Figure: runtime}
\end{figure}

In Figure~\ref{Figure: Diagonal} we also 
present the values of the empirical assortativity coefficients 
$\hat \rho_q$. 
Note that they slightly deviate from the actual 
assortativity coefficients because of the 
randomness in the construction and the erasure procedure. 
Nevertheless, by Theorem~\ref{Theorem: simple} and by 
the continuous mapping theorem, the empirical 
assortativity coefficient $\hat \rho_q$ converges in probability 
to $\rho_q$ as $N\to \infty$. 
In Figure~\ref{Figure: runtime} we analyze how much $\hat \rho_q$ 
deviates from $\rho_q$ for the fixed distributions $P_{0.5}$ and $Q_{q}$ from above 
with $q=14/15$, so that $\rho_q=0.8$. 
For different values of $N$, we simulate $100$ (simple) graphs of size $N$ 
using the implementation of the algorithm provided in~\eqref{Link} with $\delta=0.5+0.0001$.
Figure~\ref{Figure: runtime} illustrates that the convergence 
of $\hat \rho_q$ to the true value $\rho_q$ is reasonably fast, while  
we observe a bias for smaller values of $N$. 
Finally, Figure~\ref{Figure: runtime} indicates that the runtime 
of one graph simulation is approximately linear in $N$ 
for the given distributions $P_{0.5}$ and $Q_{14/15}$.

\section{Proofs}\label{Section: Proofs}

We start with the proof of Lemma~\ref{Lemma: acceptance set} 
which states that $\p[A_N] \to 1$ as $N\to\infty$. 

~

\begin{Proof}[of Lemma~\ref{Lemma: acceptance set}]
For each $k\in [K]_1$, $n_k^+$ has a binomial 
distribution with parameters $N'$ and $p_k^+$. 
Therefore, by Chebyshev's inequality,  
\begin{equation*}
\p\left[ \left|  n_k^+ - p_k^+ N' \right| > p_k^+N^\delta/2  \right] 
~\le~
\frac{N' p_k^+(1-p_k^+)}{N^{2\delta}(p_k^+)^2/4}
~\to~0,
\quad \text{as $N\to\infty$,} 
\end{equation*}
since $\delta \in (1/2,1)$. 
Similarly,  $\sum_{e=1}^{\lceil zN'' \rceil} 1_{\{k_e=k\}}$ has a binomial 
distribution with parameters $\lceil zN'' \rceil$ and $q_k^+$. 
Therefore, by condition~\eqref{Condition 2}, $\sum_{e=1}^{\lceil zN'' \rceil}1_{\{k_e=k\}} / k$ has mean 
$q_k^+\lceil zN'' \rceil / k = p_k^+ \lceil zN'' \rceil / z$ 
and its variance is of order $N''$. 
By Chebyshev's inequality it follows that 
\begin{equation*}
\p\left[ \left|  e_k^+ - p_k^+ N'' \right| > p_k^+N^\delta/2  \right] 
~\to~0,
\quad \text{as $N\to\infty$.} 
\end{equation*}
Similarly for $n_j^-$ and $e_j^-$, $j\in[J]_1$, 
using condition~\eqref{Condition 3}.  
\end{Proof}

~

\begin{Proof}[of Theorem~\ref{Lemma: random degrees}] 
Choose $s\in\N$ fixed and 
let $u: (\N_0\times\N_0)^{s} \to [-H,H]$ be a function bounded by $H>0$. 
For $N\in \N$, denote by $(j_1,k_{1}),\ldots,(j_N,k_{N})$ the 
node-types generated by the algorithm in Section~\ref{Section: Algorithm}. 
Define $\big( (\tilde j_1,\tilde k_{1}),\ldots,(\tilde j_N,\tilde k_{N}) \big )$,
where $(\tilde j_v,\tilde k_{v})= ( j_v, k_v)$ for all $v=1,\ldots, N'$, 
and $(\tilde j_v,\tilde k_{v})$, $v=N'+1,\ldots,N$, are independent 
random variables each having distribution~$P$. 
By the triangle inequality we have 
\begin{eqnarray*}
&&\hspace{-1cm}
\left| \E\left[u((j_{v_1},k_{v_1}),\ldots,(j_{v_s},k_{v_s}))\right] 
- \E\left[u((j'_1,k'_{1}),\ldots,(j'_s,k'_{s}))\right] \right|
\\&&\le~
\left| \E\left[u((j_{v_1},k_{v_1}),\ldots,(j_{v_s},k_{v_s}))
- u((\tilde j_{v_1}, \tilde k_{v_1}),\ldots,(\tilde j_{v_s},\tilde k_{v_s}))\right] \right|
\\&&\hspace{1cm}+~
\left| \E\left[u((\tilde j_{v_1}, \tilde k_{v_1}),\ldots,(\tilde j_{v_s},\tilde k_{v_s}))
- u((j'_1,k'_{1}),\ldots,(j'_s,k'_{s}))\right] \right|.
\end{eqnarray*}
Since on $A_N$, 
$\big( (\tilde j_{v_1},\tilde k_{v_1}),\ldots,(\tilde j_{v_s},\tilde k_{v_s}) \big )$ has the same 
distribution as $\big( ( j'_1, k'_{1}),\ldots,( j'_s, k'_{s}) \big )$, 
the second term on the right-hand side satisfies 
\begin{eqnarray*}
&&\hspace{-1cm}
\left| \E\left[u((\tilde j_{v_1}, \tilde k_{v_1}),\ldots,(\tilde j_{v_s},\tilde k_{v_s}))
- u((j'_1,k'_{1}),\ldots,(j'_s,k'_{s}))\right] \right|
\\&&\le~
\left| \E\left[\left. u((\tilde j_{v_1}, \tilde k_{v_1}),\ldots,(\tilde j_{v_s},\tilde k_{v_s}))
- u((j'_1,k'_{1}),\ldots,(j'_s,k'_{s}))\right| A_N \right] \right| 
	+ 2H \p[A_N^c]
\\&&=~
2H \p[A_N^c],
\end{eqnarray*}
which converges to $0$ as $N\to\infty$ by Lemma~\ref{Lemma: acceptance set}. 
For the first term we have by the definition of the random variables $(\tilde j_v,\tilde k_v)$, 
$v=1,\ldots,N$,  
\begin{eqnarray*}
&&\hspace{-1cm}
\left| \E\left[u((j_{v_1},k_{v_1}),\ldots,(j_{v_s},k_{v_s}))
- u((\tilde j_{v_1}, \tilde k_{v_1}),\ldots,(\tilde j_{v_s},\tilde k_{v_s}))\right] \right|
\\&&=~
\left| \E\left[\left(u((j_{v_1},k_{v_1}),\ldots,(j_{v_s},k_{v_s}))
- u((\tilde j_{v_1}, \tilde k_{v_1}),\ldots,(\tilde j_{v_s},\tilde k_{v_s}))\right)
1_{\bigcup_{l=1}^s\left\{v_l \in\{N'+1,\ldots,N\}   \right\}}\right] \right|
\\&&\le~
2H\sum_{l=1}^s \p\left[ v_l \in\{N'+1,\ldots,N\}   \right]
~=~
2Hs\frac{N-N'}{N},
\end{eqnarray*}
which converges to $0$ as $N\to\infty$ by the choice of $N'$. 
The corresponding result for the edge-types follows by  
exactly the same arguments since, on event $A_N$,  
the number of generated edge-types having distribution $Q$ 
is $\lceil zN'' \rceil$ and 
the number of artificially added  
edge-types is at most $K^2+J^2$, see Step~2 of the algorithm. 
\end{Proof}

~

\begin{Proof}[of Theorem~\ref{Lemma: empirical distributions 2}]
Let $j\in[J]_0$, $k\in[K]_0$ and choose $\varepsilon>0$.  
For $N$ so large that $(N-N')/N \le \varepsilon /2$ 
we have  
\begin{eqnarray*}
\p\left[ \left|
\frac{1}{N}\sum_{v=1}^N 1_{\{ j_v=j,k_v=k \}} - p_{j,k}
\right| > \varepsilon   \right] 
&\le&
\p\left[\left.  \frac{N-N'}{N} +  \left|
\frac{1}{N}\sum_{v =1}^{N'} 1_{\{ j_v=j,k_v=k \}} - p_{j,k}
\right| > \varepsilon \right| A_N  \right] 
+
\p\left[A_N^c  \right] 
\\&\le&
\p\left[\left.   \left|
\sum_{v =1}^{N'} 1_{\{ j_v=j,k_v=k \}} - p_{j,k}N
\right| > N\varepsilon/2 \right| A_N  \right] 
+ 
\p\left[A_N^c  \right]. 
\end{eqnarray*}
By Lemma~\ref{Lemma: acceptance set} it remains to consider 
the first term on the right-hand side. 
By the triangle and Chebyshev's inequality it follows that 
\begin{eqnarray*}
\p\left[\left. \left|
\sum_{v=1}^{N'} 1_{\{ j_v=j,k_v=k \}} - p_{j,k}N
\right| > N\varepsilon/2 \right| A_N  \right] 
&\le&
\frac{N' p_{j,k}(1-p_{j,k})}
{\p\left[ A_N \right]N^2 \varepsilon^2/16}
+ 
\p\left[\left.   
(N-N')
p_{j,k} > N\varepsilon/4 \right| A_N  \right],
\end{eqnarray*}
which converges to $0$ as $N\to\infty$ 
by Lemma~\ref{Lemma: acceptance set}. 
Similarly for the edge-types. 
\end{Proof}

~

In order to prove Theorem~\ref{Theorem: simple}, 
we first show that the expected number 
of self-loops and multiple edges arising from the 
construction in Section~\ref{Section: Algorithm} is bounded in $N$. 

~

\begin{lemma}\label{Lemma: expected}
Let $S_N$ be the number of self-loops and $M_N$ be the 
number of multiple edges of the graph generated by the algorithm 
in Section~\ref{Section: Algorithm}. 
There exists a finite constant $C>0$ such that 
\begin{equation*}
\E\left[\left. S_N +M_N \right| A_N\right] ~\le~ C.  
\end{equation*}
\end{lemma}

~

\begin{Proof}[of Lemma~\ref{Lemma: expected}] 
Let $v\in\{1,\ldots,N\}$ and denote by $s_v$ the number of edges $e$ with $e=\langle v,v\rangle$.  
Note that in order to have a self-loop $e=\langle v,v\rangle$, both ends of an edge 
of type $(k_v,j_v)$ have to be assigned to node $v$, see also Step~4 of the algorithm. 
Therefore, since there are $k_v$ edges leaving from $v$, we obtain 
upper bound  
\begin{eqnarray*}
\E[\left.S_N\right| A_N] 
&=& 
\sum_{v=1}^N 
\E[\left.s_v\right| A_N] 
~=~ 
\sum_{v=1}^N 
\E\left[\left. k_v\frac{\mcE_{k_v,j_v} }{\mcE_{k_v}^+}\frac{1}{\mcV_{j_v}^-}   \right| A_N\right] 
~\le~ 
\sum_{v=1}^N 
\E\left[\left.\frac{k_v}{\mcV_{j_v}^-}\right| A_N\right]. 
\end{eqnarray*}
Here, $\mcE_{k_v}^+$ denotes the number of generated edges $e=\langle v',w\rangle$ 
with out-degree of $v'$ being $k_v$ and $\mcV_{j_v}^-$ denotes the number of nodes 
having in-degree $j_v$. 
On event $A_N$, the number of nodes having 
in-degree $j_v$ is at least $\max\{1,p_{j_v}^-N'-p_{j_v}^-N^\delta/2\}$. 
It follows that 
\begin{eqnarray*}
\E[\left.S_N\right| A_N] 
&\le& 
\sum_{v=1}^N 
\E\left[\left.\frac{k_v}
{ \max\left\{1,p_{j}^-N'-p_{j}^-N^\delta/2 \right\}}
\right| A_N\right]
~\le~
\frac{N K }
{ \max\left\{1,\min_{j\in[J]_1}\left\{p_{j}^-N'-p_{j}^-N^\delta/2\right\} \right\}}. 
\end{eqnarray*}
Since $z>0$, there exists $j\in[J]_1$ with $p_j^- > 0$ and, hence, the right-hand side 
is bounded in $N$.

~

To bound  the expectation of $M_N$, 
let $v\in\{1,\ldots,N\}$ and denote by $m_v$ the number of multiple edges 
leaving from node $v$. 
The probability that two distinct edges leaving from $v$ are arriving at the same 
node $w\neq v$ is at most 
\begin{equation*}
\frac{j_w(j_w-1)}{j_w\mcV_{j_w}^-\left(j_w\mcV_{j_w}^--1\right)}1_{\{ j_w\ge2 \}}.
\end{equation*}
It follows that 
\begin{eqnarray*}
\E[\left.M_N\right| A_N] 
&=& 
\sum_{v=1}^N
\E[\left.m_v\right| A_N] 
~\le~
\sum_{v=1}^N
\sum_{w=1}^N
\E\left[\left. 
\binom{k_v}{2}
\frac{j_w(j_w-1)}{j_w\mcV_{j_w}^-\left(j_w\mcV_{j_w}^--1\right)} 
1_{\{ j_w\ge2 \}}  
1_{\{ k_v\ge2 \}}  
\right| A_N\right] 
\\&\le&
\sum_{v=1}^N
\sum_{w=1}^N
\E\left[\left. 
\frac{k_v^2 j_w^2}{2j_w\mcV_{j_w}^-\left(j_w\mcV_{j_w}^--1\right)}
1_{\{ j_w\ge2 \}} 
1_{\{ k_v\ge2 \}}  
\right| A_N\right]. 
\end{eqnarray*}
Since $2j_w\mcV_{j_w}^-\left(j_w\mcV_{j_w}^--1\right) \ge \left(j_w\mcV_{j_w}^-\right)^2$ 
for $j_w \ge2$, it follows that 
\begin{eqnarray*}
\E[\left.M_N\right| A_N] 
&\le& 
\sum_{v=1}^N
\sum_{w=1}^N
\E\left[\left. 
\frac{k_v^2}{\left(\mcV_{j_w}^-\right)^2}1_{\{ j_w\ge2 \}} 1_{\{ k_v\ge2 \}}    \right| A_N\right]. 
\end{eqnarray*}
Using that on $A_N$, the number of nodes having 
in-degree $j_w$ is at least $\max\{1,p_{j_w}^-N'-p_{j_w}^-N^\delta/2\}$, 
it follows that 
\begin{eqnarray*}
\E[\left.M_N\right| A_N] 
&\le& 
\frac{N^2K^2}
{\max\left\{1,\min_{j\in[J]_1}\left\{p_{j}^-N'-p_{j}^-N^\delta/2\right\} \right\}^2}. 
\end{eqnarray*}
The right-hand side is again bounded in $N$. 
This finishes the proof of Lemma~\ref{Lemma: expected}. 
\end{Proof}

~

\begin{Proof}[of Theorem~\ref{Theorem: simple}]
We first show that Theorem~\ref{Lemma: empirical distributions 2} 
holds true for the erased configuration graph. 
Let $S_N$ be the number of self-loops and  
let $M_N$ be the number of multiple edges generated by the algorithm. 
For $j \in [J]_0$ and $k\in [K]_0$ denote by $\mcV_{j,k}^e$ 
the number of constructed nodes of type $(j,k)$ after 
erasing all self-loops and multiple edges. 
Similarly we define $\mcE_{k,j}^e$ for $k \in [K]_1$ and $j\in [J]_1$. 
In order to prove that Theorem~\ref{Lemma: empirical distributions 2} 
holds true for the erased configuration graph, 
we show that for any $\varepsilon>0$, 
\begin{equation}\label{Equation: erased empirical}
\lim_{N\to\infty}
\p\left[
\sum_{j \in [J]_0, k\in [K]_0}
\left|
\frac{\mcV^e_{j,k}}{N} 
-
p_{j,k}
\right| 
\ge \varepsilon
\quad \text{ or } \quad
\sum_{k \in [K]_1, j\in [J]_1}
\left|
\mcE^e_{k,j}
-
q_{k,j}\mcE^e
\right|
\ge \varepsilon \mcE^e
\right]
~=~0,
\end{equation}
where $\mcE^e$ denotes the total number of edges 
in the erased configuration graph (which could be equal to $0$ if all edges 
of the constructed graph are self-loops).  
Let $\varepsilon >0$ and choose $j \in [J]_0$ and $k\in [K]_0$. 
We have 
\begin{equation*}
\p\left[
\left|
\frac{\mcV^e_{j,k}}{N} 
-
p_{j,k}
\right|
\ge \varepsilon
\right]
~\le~
\p\left[
\left|
\frac{\mcV^e_{j,k}- \mcV_{j,k}}{N} 
\right|
\ge \varepsilon/2
\right]
+
\p\left[
\left|
\frac{\mcV_{j,k}}{N} 
-
p_{j,k}
\right|
\ge \varepsilon/2
\right].
\end{equation*}
The second term on the right-hand side converges to $0$ as $N\to\infty$ 
by Theorem~\ref{Lemma: empirical distributions 2}. 
For the first term, note that 
\begin{equation*}
\p\left[
\left|
\frac{\mcV^e_{j,k}- \mcV_{j,k}}{N} 
\right|
\ge \varepsilon/2
\right]
~\le~
\p\left[
\sum_{v=1}^N 
\left|1_{\{j_v^e=j,k_v^e=k\}}- 1_{\{j_v=j,k_v=k\}}
\right|
\ge \varepsilon N/2
\right],
\end{equation*}
where $(j_v^e,k^e_v)$ denotes the type of node $v$ in the erased configuration graph. 
Denote by $s_v$ the number of self-loops attached to $v$.  
Denote by $m_v^+$ the number of multiple edges leaving from $v$ and denote 
by $m_v^-$ the number of multiple edges arriving at $v$. 
Note that $j_v^e=j_v$ if and only if $s_v+m^-_v=0$, and similarly 
$k_v^e=k_v$ if and only if $s_v+m^+_v=0$.  
We therefore have that 
\begin{eqnarray*}
1_{\{j_v^e=j,k_v^e=k\}}- 1_{\{j_v=j,k_v=k\}}
&=&
1_{\{ s_v+m_v^->0 \text{ or } s_v + m_v^+>0 \}}\left(1_{\{j_v^e=j,k_v^e=k\}}- 1_{\{j_v=j,k_v=k\}}\right)
\\&=&
1_{\{ s_v+m_v>0\}}\left(1_{\{j_v^e=j,k_v^e=k\}}- 1_{\{j_v=j,k_v=k\}}\right),
\end{eqnarray*}
where we set $m_v=m_v^+ + m_v^-$. 
Hence, 
\begin{eqnarray*}
\p\left[
\left|
\frac{\mcV^e_{j,k}- \mcV_{j,k}}{N} 
\right|
\ge \varepsilon/2
\right]
&\le&
\p\left[
\sum_{v=1}^N 
1_{\{ s_v+m_v>0\}}
\ge \varepsilon N/2
\right]
~\le~
\p\left[
\sum_{v=1}^N 
(s_v+m_v)
\ge \varepsilon N/2
\right].
\end{eqnarray*}
It follows by Markov's inequality, Lemma~\ref{Lemma: acceptance set} and 
Lemma~\ref{Lemma: expected} that 
\begin{eqnarray*}
\p\left[
\left|
\frac{\mcV^e_{j,k}- \mcV_{j,k}}{N} 
\right|
\ge \varepsilon/2
\right]
&\le&
\p\left[
S_N+2M_N
\ge N\varepsilon/2
\right]
~\le~
\frac{\E\left[\left. S_N+2M_N \right| A_N \right]}{N\varepsilon/2} 
+
\p\left[A_N^c\right]
~\to~0,
\quad \text{ as $N\to\infty$.}
\end{eqnarray*}
We now prove the same result for the edge-types 
under the conditional probability, conditional given $A_N$, 
which is enough due to Lemma~\ref{Lemma: acceptance set}. 
Note that on event $A_N$ the number of generated edges $\mcE$ is at least 
$zN''$, see Step~1 of the algorithm. 
It follows by Markov's inequality and Lemma~\ref{Lemma: expected}, 
and since $\mcE - \mcE^e =  S_N + M_N$, 
that for every $\eta>0$, 
\begin{equation}\label{Equation: mcE}
\p\left[\left.
\left|\frac{\mcE^e}{\mcE} -1 \right|
\ge \eta
\;\right| A_N
\right]
~=~
\p\left[\left.
\frac{S_N+M_N}{\mcE} 
\ge \eta
\;\right| A_N
\right]
~\le~
\p\left[\left.
S_N+M_N
\ge zN''\eta
\;\right| A_N
\right]
~\to~0,
\end{equation}
as $N\to \infty$. 
Note that 
for fixed $k \in [K]_1$ and $j\in [J]_1$, 
\begin{eqnarray*}
\p\left[\left.
\left|
\mcE^e_{k,j}
-
q_{k,j}\mcE^e
\right|
\ge \varepsilon\mcE^e
\;\right| A_N
\right]
&\le&
\p\left[\left.
\left|
\mcE^e_{k,j}- \mcE_{k,j}
\right|
\ge \varepsilon\mcE^e/2
\;\right| A_N
\right]
+
\p\left[\left.
\left|
\mcE_{k,j}
-
q_{k,j}\mcE^e
\right|
\ge \varepsilon\mcE^e/2
\;\right| A_N
\right].
\end{eqnarray*}
For the second term on the right-hand side we have
\begin{eqnarray*}
\p\left[\left.
\left|
\mcE_{k,j}
-
q_{k,j}\mcE^e
\right|
\ge \varepsilon\mcE^e/2
\;\right| A_N
\right]
&=&
\p\left[\left.
\left|
\frac{\mcE_{k,j}}{\mcE} 
-
\frac{\mcE - S_N - M_N}{\mcE} q_{k,j}
\right|
\ge \frac{\mcE^e}{\mcE}\varepsilon/2
\;\right| A_N
\right]
\\&&\hspace{-3cm}\le~
\p\left[\left.
\left|
\frac{\mcE_{k,j}}{\mcE} 
-
q_{k,j}
\right|
\ge \frac{\mcE^e}{\mcE}\varepsilon/4
\;\right| A_N
\right]
+
\p\left[\left.
\left(S_N +M_N\right) q_{k,j}
\ge \mcE^e\varepsilon/4
\;\right| A_N
\right]
\\&&\hspace{-3cm}\le~
\p\left[\left.
\left|
\frac{\mcE_{k,j}}{\mcE} 
-
q_{k,j}
\right|
\ge \frac{\mcE^e}{\mcE}\varepsilon/4
\;\right| A_N
\right]
+
\p\left[\left.
S_N +M_N
\ge \frac{zN''\varepsilon/4}{q_{k,j}+ \varepsilon/4}
\;\right| A_N
\right],
\end{eqnarray*}
where in the last step we used that $\mcE^e =  \mcE - (S_N + M_N) \ge  zN''  - (S_N + M_N)$. 
By Theorem~\ref{Lemma: empirical distributions 2}, 
\eqref{Equation: mcE} and Lemma~\ref{Lemma: expected}, 
the right-hand side converges to $0$ as $N\to\infty$.  
To bound the first term, note that the erasure procedure changes 
the number of edges of type $(k,j)$ because such edges may get erased, 
but also because an erased edge $e=\langle v,w \rangle$ changes the types of 
unerased edges that are leaving from node $v$ or that are arriving at node $w$. 
It follows that 
\begin{equation*}
\p\left[\left.
\left|
\mcE^e_{k,j}- \mcE_{k,j}
\right|
\ge \varepsilon\mcE^e/2
\;\right| A_N
\right]
~\le~
\p\left[\left.
S_N+M_N + 
\sum_{e=1}^{\mcE^e}
\left|
1_{\{ k^e_e=k, j_e^e=j \}}- 1_{\{ k_e=k, j_e=j \}}
\right|
\ge \varepsilon\mcE^e/2
\;\right| A_N
\right],
\end{equation*}
where $(k^e_e, j_e^e)$ denotes the type of edge $e$  
in the erased configuration graph. 
For an edge $e=\langle v,w \rangle$ we have $k_e^e = k_e$ if and only if 
$s_v + m_v^+=0$, and similarly $j_e^e = j_e$ if and only if 
$s_w + m_w^-=0$. 
Therefore, for $e=\langle v,w \rangle$, 
\begin{eqnarray*}
1_{\{ k^e_e=k, j_e^e=j \}}- 1_{\{ k_e=k, j_e=j \}}
&=&
1_{\{ s_v+m_v^+>0 \text{ or } s_w + m_w^->0 \}}
	\left(1_{\{ k^e_e=k, j_e^e=j \}}- 1_{\{ k_e=k, j_e=j \}}\right)
\\&=& 
1_{\{ s_e+m_e>0\}}
	\left(1_{\{ k^e_e=k, j_e^e=j \}}- 1_{\{ k_e=k, j_e=j \}}\right),
\end{eqnarray*}
where we set $s_e = s_v+s_w$ and $m_e=m_v^++m_w^-$. 
It follows that 
\begin{eqnarray*}
\p\left[\left.
\left|
\mcE^e_{k,j}- \mcE_{k,j}
\right|
\ge \varepsilon\mcE^e/2
\;\right| A_N
\right]
&\le&
\p\left[\left.
S_N+M_N+ 
\sum_{e=1}^{\mcE^e}
1_{\{ s_e+m_e>0\}}
\ge \varepsilon\mcE^e/2
\;\right| A_N
\right]
\\&&\hspace{-3cm}\le~
\p\left[\left.
S_N+M_N + 
\sum_{v=1}^N\left( k_v1_{\{ s_v+m^+_v>0\}} + j_v1_{\{ s_v+m_v^->0\}}  \right)
\ge \mcE^e \varepsilon/2
\;\right| A_N
\right]
\\&&\hspace{-3cm}\le~
\p\left[\left.
S_N+M_N
+
K (S_N+M_N) 
+ J (S_N+M_N) 
\ge \mcE^e\varepsilon/2
\;\right| A_N
\right],
\end{eqnarray*}
which converges to $0$ as $N\to \infty$ by Lemma~\ref{Lemma: expected} 
and the fact that $\mcE^e =  \mcE - (S_N + M_N) \ge  zN'' - (S_N + M_N)$. 
This finally proves~\eqref{Equation: erased empirical}. 

~

We now prove that Theorem~\ref{Lemma: random degrees}  holds 
true for the erased configuration graph. 
Denote by $\mcV_{\rm mod}$ the set of nodes whose types have  
been changed due to the erasure procedure. 
For the probability that a uniformly chosen node $v$ 
belongs to $\{N'+1,\ldots,N\}\cup\mcV_{\rm mod}$ we have 
\begin{eqnarray*}
\p\left[  v \in \{N'+1,\ldots,N\}\cup\mcV_{\rm mod}  \right] 
&=&
\E\left[  \frac{ \left| \{N'+1,\ldots,N\}\cup\mcV_{\rm mod}\right|}{N}   \right] 
~\le~
\E\left[  
\frac{ N-N' +   S_N+2M_N }
	{N}  \right]. 
\end{eqnarray*}
Since the graph returned by the algorithm in case event $A_N$ 
does not hold has no self-loops 
or multiple edges, see Step~5, it follows that 
\begin{equation*}
\p\left[  v \in \{N'+1,\ldots,N\}\cup\mcV_{\rm mod}  \right] 
~\le~
\E\left[\left.  
\frac{ N-N' +   S_N+2M_N }
	{N}  \right| A_N \right] 
+
\frac{ N-N' }{N},		
\end{equation*}
which converges to $0$ as $N\to\infty$ by  Lemma~\ref{Lemma: expected}. 
Going through the proof of Theorem~\ref{Lemma: random degrees}, 
we see that this observation is enough to conclude 
that the types of $s\in\N$ randomly chosen nodes of  
the erased configuration graph converge in distribution 
to a sequence of $s$ independent random variables each having distribution $P$ as $N\to\infty$. 
To prove the same result for the edge-types 
of the erased configuration graph, 
note that it may happen that all edges of the graph generated by the 
algorithm of Section~\ref{Section: Algorithm} are self-loops, i.e.~$\mcE^e=0$. 
In this case the edge set is 
empty for the erased configuration graph and 
we define ``the type of a randomly 
chosen edge'' to be identical to $(1,1)$ if event $\{\mcE^e=0\}$ occurs, and 
we define it to be as usual if event $\{\mcE^e=0\}$ does not occur. 
Nevertheless, the probability of event $\{\mcE^e=0\}$ 
converges to $0$ as $N\to\infty$ by Lemma~\ref{Lemma: acceptance set} 
and~\eqref{Equation: mcE} above. 
Therefore, using similar arguments as above for the node-types, we 
conclude that the types of $s\in\N$ randomly chosen edges of  
the erased configuration graph converge in distribution 
to a sequence of $s$ independent random variables each having 
distribution $Q$ as $N\to\infty$.
\end{Proof}

\begin{appendix}
\renewcommand{\thesection}{A}
\section{Variations of edge-types}\label{Appendix}

We have defined assortativity through edge-types, which, 
for an edge $e=\langle v,w \rangle$ connecting 
node $v$ to node $w$, is defined to be the tuple  
$(k_e,j_e)$ with $k_e$ denoting the out-degree of  $v$ and 
$j_e$ denoting the in-degree of  $w$. 
There are three other possibilities to define the type of an edge $e$, for instance, 
by the tuple $(k_e,k'_e)$ with $k_e$ and $k'_e$ denoting the out-degree of  $v$ and 
$w$, respectively. 
We briefly explain how we need to modify the presented algorithm in Section~\ref{Section: Algorithm} 
when redefining the type of an edge as above.  
In the same spirit one can then construct algorithms for the 
remaining two definitions of edge-types. 
Results as in Lemma~\ref{Lemma: acceptance set} and 
Theorems~\ref{Lemma: random degrees}--\ref{Theorem: simple} 
can also  be proven using the new definition of edge-types.

~

Define the type of an edge $e=\langle v,w \rangle$ 
by the tuple $(k_e,k'_e)$ with $k_e$ and $k'_e$ denoting the out-degree of node $v$ and 
$w$, respectively. 
In view of conditions~\eqref{Condition 2} and~\eqref{Condition 3} the corresponding 
edge-type distribution $Q=(q_{k,k'})$ with marginal distributions  
$(q_k^l)_{k\in[K]_1}$ and $(q_{k'}^r)_{k'\in[K]_0}$ then needs to satisfy, 
for given node-type distribution $P$ with mean degree $z>0$, 
\begin{equation}\label{Equation: variation}
q_k^l~=~kp_k^+/z, \quad k\in[K]_1,
\qquad \text{and} \qquad 
q_{k'}^r~=~\sum_{j\in[J]_1}jp_{j,k'}/z, \quad k'\in[K]_0.
\end{equation}
Here, superscript ``$l$'' refers to node $v$ of the edge $e=\langle v,w \rangle$ and 
superscript ``$r$''  to node $w$. 
Condition~\eqref{Equation: variation} is justified by counting nodes and 
edges of corresponding types in a given graph. 
For instance, note that for every node of type $(j,k')$,  $j\in [J]_1$ and 
$k'\in [K]_0$, we have $j$ edges $e=\langle v,w \rangle$ 
with out-degree of $w$ being $k'$. 

~

We now describe a modification of the algorithm in Section~\ref{Section: Algorithm} 
to construct graphs corresponding to the new distribution $Q$. 
Choose $N\in\N$ so large that there 
exists  $N'\in\N$ with  
$N=N' + 2(1+z)\left\lceil N^\delta \right \rceil + K^2$,  
and set $N''=N' + \left\lceil N^\delta\right \rceil$. 
\begin{enumerate}[label={}, leftmargin=10pt, rightmargin=10pt]
\item 
{\bf Step 1.}  
Assign to each node $v=1,\ldots,N'$ independently a node-type $(j_v,k_v)$ 
according to distribution $P$. 
Generate edges $e=1,\ldots, \lceil zN'' \rceil$ 
having independent edge-types $(k_e,k'_e)$ according to distribution $Q$,  
independently of the node-types. 
Define
\begin{eqnarray*}
n_k^l &=& \sum_{v=1}^{N'} 1_ {\{ k_v=k \}}
\quad \text{ and } \quad 
e_k^l ~=~ \left\lceil \frac{1}{k}\sum_{e=1}^{\lceil zN'' \rceil} 1_ {\{ k_e=k \}}\right\rceil
\quad \text{ for all $k\in[K]_1$};
\\
n_{k'}^r &=& \sum_{j=1}^{J}j\sum_{v=1}^{N'} 1_ {\{ j_v=j,\,k_v=k' \}}
\quad \text{ and } \quad 
e_{k'}^r ~=~ \sum_{e=1}^{\lceil zN'' \rceil} 1_ {\{ k'_e=k' \}}
\quad\text{ for all $k'\in[K]_0$}.  
\end{eqnarray*}
Let $A_N$ be the event on which we have, 
set $p_{k'}^r=\sum_{j=1}^{J}jp_{j,k'}$, 
\begin{eqnarray*}
\left|  n_k^l - p_k^+ N' \right| \le p_k^+N^\delta/2
&\quad \text{ and } \quad &
\left|  e_k^l - p_k^+ N'' \right| \le p_k^+N^\delta/2
\quad \text{ for all $k\in[K]_1$}; 
\\
\left|  n_{k'}^r - p_{k'}^r N' \right| \le p_{k'}^rN^\delta/2
&\quad \text{ and } \quad &
\left|  e_{k'}^r - p_{k'}^r N'' \right| \le p_{k'}^rN^\delta/2
\quad\text{ for all $k'\in[K]_0$}.  
\end{eqnarray*}
Only accept Step~1 if event $A_N$ occurs and proceed to Step 2, otherwise return a 
graph containing at least one edge. 
Observe that by the relation between $N'$ and $N''$ we have, on $A_N$, $n_k^l - e_k^l \le 0$ 
for all $k\in[K]_1$ and $n_{k'}^r - e_{k'}^r \le 0$ for all $k'\in[K]_0$.

\item
{\bf Steps 2 and 3.}
If the number of generated 
edges $e=\langle v,w \rangle$ with out-degree of $v$ 
being $k\in[K]_1$ is not a multiple of $k$, we add the corresponding number of 
edges of type $(k,0)$ and nodes of type $(1,0)$. 
Additionally, we add $e_k^l - n_k^l \ge 0$ nodes 
of type $(0,k)$, i.e.~we finally obtain that the number of nodes 
with out-degree $k$ corresponds to the number of such nodes needed for the  
edges $e=\langle v,w \rangle$ with out-degree of $v$ 
being $k$. 

Note that $n_{k'}^r$ is the sum of in-degrees 
of all nodes generated in Step~1 having out-degree $k'\in[K]_0$, 
which is, since $e_{k'}^r - n_{k'}^r \ge 0$, 
eventually strictly less than $e_{k'}^r$. 
Therefore, we add $e_{k'}^r - n_{k'}^r$ nodes of type $(1,k')$, 
and for each such node we add an additional node of type $(k',0)$.
The node of type $(k',0)$ is only needed if $k'\ge1$, because 
for each added node of type $(1,k')$ we additionally 
add $k'$ edges of type $(k',0)$. 
Note that instead of adding a node of type $(k',0)$ we could also add $k'$ nodes of type $(1,0)$. 

\item
{\bf Step 4.}
For each $k\in[K]_1$, assign to each node 
having out-degree $k$ exactly $k$ uniformly chosen  
edges $e$ of type $(k_e,k'_e)$ with $k_e=k$. 
For each $k'\in[K]_0$, assign to every node $v$ of type $(j_v,k_v)$ with $k_v=k'$ 
exactly $j_v$ uniformly chosen 
edges $e$ of type $(k_e,k'_e)$ with $k'_e=k'$. 
Return the constructed graph. 
\end{enumerate}

~
 
 \noindent
Following this modified algorithm, the total number of nodes we 
add in Steps~2 and~3 is at most, 
on event $A_N$, 
\begin{equation*}
\sum_{k=1}^K  \left(k-1 + (e_k^l - n_k^l) \right) 
+ 
\sum_{k'=0}^K  (e_{k'}^r - n_{k'}^r)\left( 1+1 \right)
~\le~
K^2 + 2\lceil N^\delta  \rceil + 2z\lceil N^\delta  \rceil. 
\end{equation*}
This implies that we have sufficiently many undetermined nodes $N-N'$ to which  
we can assign corresponding node-types. 
Moreover, the total number of edges we add in Steps~2 and~3 is at most, 
on event $A_N$, 
\begin{equation*}
\sum_{k=1}^K  \left(k-1\right) 
+ 
\sum_{k'=0}^K  k'(e_{k'}^r - n_{k'}^r) 
~\le~
K^2 + zK \lceil N^\delta  \rceil. 
\end{equation*}

~

\begin{remark}
Note that when considering a variation of edge-types with
corresponding distribution $Q$ 
satisfying~\eqref{Equation: variation} or variations thereof, 
then the marginal distributions of $Q$ are determined by the distribution $P$.  
Therefore, the discussion in the beginning of Section~\ref{Section: Examples} 
carries over also for variations of edge-types. 
\end{remark}

\end{appendix}


\end{document}